\documentclass[10pt,reqno]{amsart}

\usepackage{amssymb,amsmath,mathrsfs}

\newif\iffigures\figurestrue
\figuresfalse
\newif\ifhyper\IfFileExists{hyperref.sty}{\hypertrue}{\hyperfalse}
\ifhyper\usepackage{hyperref}
\def\hitem#1#2{\item[\hypertarget{#1}{#2}]\expandafter\gdef\csname LBL#1ITM\endcsname{#2}}
\def\iref#1{\hyperlink{#1}{\csname LBL#1ITM\endcsname}}
\else
\def\hitem#1#2{\item[{#2}]\expandafter\gdef\csname LBL#1ITM\endcsname{#2}}
\def\iref#1{{\csname LBL#1ITM\endcsname}}
\fi

\def\bl{\begin{lemma}}
\def\el{\end{lemma}}
\def\bth{\begin{theorem}}
\def\eth{\end{theorem}}
\def\bc{\begin{corollary}}
\def\ec{\end{corollary}}
\def\bcj{\begin{conjecture}}
\def\ecj{\end{conjecture}}
\def\bpr{\begin{proposition}}
\def\epr{\end{proposition}}
\def\bde{\begin{definition}}
\def\ede{\end{definition}}
\def\bq{\begin{question}}
\def\eq{\end{question}}

\newcommand{\be}{\begin{eqnarray}}
\newcommand{\ee}{\end{eqnarray}}

\def\Pr{\mathbb{P}}

\def\st{{\, : \ }}
\def\gen#1{\langle#1\rangle}
\def\Gen#1{\bigl\langle#1\bigr\rangle}
\def\supp{\mathop{\rm supp}}
\let\phi=\varphi
\def\G{\mathsf{G}}
\def\BS{\mathsf{BS}}
\def\F{\mathsf{F}}
\def\St{\mathrm{St}}
\def\p{\partial}
\def\Aut{\mathrm{Aut}}
\def\shift{\textsf{s}}

\def\eps{\varepsilon}

\def\Z{{\mathbb Z}}
\def\N{{\mathbb N}}
\def\T{{\mathcal T}}

\def\NN{{\mathcal N}}
\def\JJ{{\mathcal J}}
\def\SS{{\mathcal S}}
\def\LL{{\mathcal L}}

\def\Cl{{\mathscr C}}

\newtheorem{theorem}{Theorem}[section]
\newtheorem{definition}{Definition}[section]
\newtheorem{lemma}[theorem]{Lemma}
\newtheorem{corollary}[theorem]{Corollary}
\newtheorem{proposition}[theorem]{Proposition}
\newtheorem{conjecture}[theorem]{Conjecture}
\newtheorem{question}{Question}[section]

\theoremstyle{definition}
\numberwithin{equation}{section}

\input epsf.sty

\begin{document}

\title{Scale-invariant groups}

\author{Volodymyr Nekrashevych}
\address{Department of Mathematics, Texas A\&M University,
College Station, TX 77843-3368, USA}
\email{nekrash@math.tamu.edu}
\urladdr{http://www.math.tamu.edu/$\sim$nekrash}

\author{G\'abor Pete}
\address{Department of Mathematics, University of Toronto, 40 St George St., Toronto, ON, M5S 2E4, Canada}
\email{gabor@math.toronto.edu; robagetep@gmail.com}
\urladdr{http://www.math.toronto.edu/$\sim$gabor}

%\date{\today}

% \keywords{finite automata, self-similar groups, expanding homomorphism, renormalization, percolation, lamplighter group, solvable Baumslag-Solitar group, affine group, Diestel-Leader graphs, co-Hopfian groups}

\maketitle

\begin{abstract} Motivated by the renormalization method in statistical physics, Itai Benjamini defined a finitely generated infinite group $G$ to be scale-invariant if there is a nested sequence of finite index subgroups $G_n$ that are all isomorphic to $G$ and whose intersection is a finite group. He conjectured that every scale-invariant group has polynomial growth, hence is virtually nilpotent. We disprove his conjecture by showing that the following groups (mostly finite-state self-similar groups) are scale-invariant: the lamplighter groups $\F\wr \Z$, where $\F$ is any finite Abelian group; the solvable Baumslag-Solitar groups $\BS(1,m)$; the affine groups $A\ltimes \Z^d$, for any $A\leq GL(\Z,d)$. However, the conjecture remains open with some natural stronger notions of scale-invariance for groups and transitive graphs. We construct scale-invariant tilings of certain Cayley graphs of the discrete Heisenberg group, whose existence is not immediate just from the scale-invariance of the group. We also note that torsion-free non-elementary hyperbolic groups are not scale-invariant.
\end{abstract}

\section{Introduction}\label{s.intro}

Itai Benjamini introduced the following notion. A finitely generated infinite group $G$ is called {\bf scale-invariant} {\it if there is a nested sequence of finite index subgroups $G_n$ that are all isomorphic to $G$ and whose intersection is a finite group}. He conjectured on his website \cite{Itai} that every scale-invariant group has polynomial growth, and hence, by \cite{Gromov}, is virtually nilpotent. (His definition was originally slightly weaker, and that was the form in which the conjecture was popularized by Mark Sapir's survey of open problems \cite[Problem 9.12]{Sapir}.)

The main motivation for defining this notion comes from statistical mechanics. A key tool in the study of percolation, the Ising model and other stochastic processes on $\Z^d$ is {\bf renormalization}, whose main geometric ingredient is that the lattice $\Z^d$ can be tiled by large boxes such that the resulting tiling is isomorphic to the original lattice. See \cite{Grimm} for background on percolation renormalization on $\Z^d$. In the past decade or so, statistical mechanics on Cayley graphs other than $\Z^d$ has been a lively research area, see \cite{BS:beyond,Ly:survey,LPbook}, but the renormalization technique has not been generalized yet. And, at least for the case of $\Z^d$, the existence of this tiling seems to be tied to the scale-invariant structure of the group:  
we get a good tiling by $2^{k+1}$-boxes from grouping $2^d$ neighboring $2^k$-boxes, which are cosets  of the subgroup $2^k\Z^d$, with adjacency defined by $2^k\Z^d\simeq \Z^d$, and finally, these boxes can exhaust the lattice as $k\to\infty$ because $\bigcap_{k\geq 0} 2^k \Z^d=\{\mathbf{0}\}$.
%the $n\times n$ boxes forming a tiling are cosets of the subgroup $n\Z^d$, the tiling graph is $n\Z^d\simeq \Z^d$, and  the boxes  exhaust the lattice as $n\to\infty$ because $\bigcap_{n\geq 1} n\Z^d=\{\mathbf{0}\}$.

On the other hand, there is a tempting geometric approach to prove the conjecture.
% On the other hand, one naturally comes immediately notice the reason for 
% Benjamini's conjecture was the following geometric intuition. 
It seems we cannot lose much of the geometric content of the problem by assuming that there is an injective endomorphism $\phi$ of $G$ such that the nested sequence $G_n:=\phi^{\circ n}(G)$ has the desired properties. Then $\phi$ seems to be almost an {\bf expanding homomorphism} (i.e., one that has a convolution power that increases all distances by a uniform factor larger than 1), with $[G:\phi(G)]<\infty$. However, such homomorphisms exist only in groups of polynomial volume growth, as a simple but important theorem due in different versions to \cite{Franks}, \cite{Farkas} and \cite{Gelb} says. As we will see, this argument is wrong, but we do not know which ``it seems'' step can be blamed (maybe both).   
% This result was an important ingredient in Gromov's proof that groups of polynomial growth are virtually nilpotent \cite{Gromov}.
% If one believes that such a scaling endomorphism should really be expanding,
% then it does not seem to hard to believe

In this note we disprove Benjamini's conjecture by giving several examples of scale-invariant groups that are not virtually nilpotent. We also introduce some stronger notions of scale-invariance for groups and transitive graphs; for these, the polynomial growth conjecture remains open.

\bth\label{t.general}
Let $H$ be a countable scale-invariant group, with a descending sequence $H=H_0>H_1>\ldots$ of finite index subgroups, each isomorphic to $H$, with trivial intersection. Let $A$ be a countable automorphism group of $H$ leaving all subgroups $H_n$ invariant. Assume that the action of $A$ is
faithful on each $H_n$, and that the semidirect products $A\ltimes H_n$ are isomorphic to $A\ltimes H$. Then $G:=A\ltimes H$ is scale-invariant; in fact, there is a required subgroup chain with $\bigcap_{n\geq 0} G_n=\{1\}$.
\eth

\bc\label{c.examples}
The following groups are scale-invariant.
\begin{itemize}
\item[(1)] The {\bf lamplighter groups} $\G=\F\wr\Z$, where $\F$ is any finite Abelian group.
\item[(2)] The solvable {\bf Baumslag-Solitar groups} $\BS(1,m)=\Gen{a,b \mid bab^{-1}=a^m}$ with $m > 1$.
\item[(3)]  The {\bf affine group} $GL(\Z,d) \ltimes \Z^d$, and its subgroups
$A\ltimes \Z^d$ for any $A\leq GL(\Z,d)$, $d > 1$.
\end{itemize}
\ec

The lamplighter and Baumslag-Solitar groups are solvable groups of exponential growth, which have served as interesting (counter)examples for many questions since their introduction \cite{KV,BS}. Note furthermore that the affine groups give examples that are even non-amenable.

In proving Theorem~\ref{t.general}, a key tool will be that $G=A\ltimes H$ acts naturally on an infinite rooted tree: the coset tree of the sequence $(H_n)_{n\geq 0}$. Moreover, most of our examples in Corollary~\ref{c.examples} are actually {\bf finite-state self-similar groups}, even though this is not at all obvious from their original definitions. (See \cite{fractal,Nekra} for background on self-similar groups.) For the lamplighter group (for $\F=\Z_2$) this was first noticed by Grigorchuk and \.Zuk \cite{GZ}, with a simpler proof in \cite{GNS}, and for general $\F$ by \cite{SiSt}; for the Baumslag-Solitar group by Bartholdi and \v{S}uni\'k \cite{solvauto}; for the full affine groups by Brunner and Sidki \cite{affine}. Although it has been well-known that a self-similar group may contain a finite index subgroup that is virtually the direct product of a few copies of itself (the typical examples are branch groups \cite{branch}, most notably Grigorchuk's groups), it does not  seem to have been observed before that some of the not virtually nilpotent examples contain finite index copies of themselves.

The subgroups $G_n$ will be the vertex stabilizers along an infinite ray in the rooted tree on which $G$ acts, so $\bigcap_{n\geq 1} G_n$ will be the stabilizer of the entire ray. Each $G_n$, the stabilizer of a vertex $v$, will be of the form $\phi_{v}(G)$, and if the parent of $v$ is $u$, then $\phi_v=\phi_u\circ \phi_{u,v}$, where $\phi_{u,v}$ is one of a finite number of injective endomorphisms of $G$ that is chosen according to ``which'' child of $u$ our $v$ is. In order to get a trivial intersection, we will need to take an ``irrational ray'' in some sense; in particular, we will prove in most of our examples that the stabilizer of a ``periodic ray'' (where the sequence of endomorphisms $\phi_{u,v}$ is periodic) is always infinite.  A consequence is that the following natural version of Benjamini's problem remains open:

\bq\label{q.single} Suppose that a group $G$ is {\bf strongly scale-invariant}: there exists a single injective endomorphism $\phi$ of $G$ such that $\bigcap_{n\geq 0}\phi^{\circ n}(G)$ is finite. 
Must $G$ have polynomial growth?
\eq

There are many strongly scale-invariant nilpotent groups: e.g., the {\bf integer Heisenberg group} 
$$
\left\{(x,y,z)=\begin{pmatrix}1 & x & y\\ 0 & 1 & z\\ 0 & 0 & 1\end{pmatrix} \st x,y,z\in \Z\right\},
\quad\text{with}\quad \phi(x,y,z)=(2x,4y,2z)\,.
$$
Even though these examples are covered by part (3) of Corollary \ref{c.examples}, e.g., the Heisenberg group equals $\Z\,_M\hspace{-0.1 cm}\ltimes \Z^2$ with $M=\begin{pmatrix}1 & 1 \\ 0 & 1 \end{pmatrix}$, our general construction does not give a strongly scale-invariant sequence for them: in that statement, whenever $G=A\ltimes \Z^d$ is nilpotent, $A$ must consist of unipotent matrices only (see, e.g., \cite{Kap:GGT}), and for this case we will show that periodic rays give non-trivial stabilizers. This is particularly interesting (and sad, maybe) because if there was a matrix $M_{d\times d}$ giving a strongly scale-invariant group $\Z\,_M\hspace{-0.1 cm} \ltimes \Z^d$, then it would be easy to combine it with a hyperbolic matrix $N_{\ell\times\ell}$ so that the resulting block-diagonal matrix $M\oplus N$ would yield a strongly scale-invariant group $\Z\,_{M\oplus N}\hspace{-0.1 cm}\ltimes \Z^{d+\ell}$ of exponential growth. 

From the percolation point of view, the following geometric-combinatorial version of scale-invari\-ance is also very natural. Recall that a graph $\Gamma=(V,E)$ is called transitive if its automorphism group acts transitively on the vertex set $V$.

\bq\label{q.tiling} A {\bf scale-invariant tiling} of a transitive graph $\Gamma$ is a decomposition of its vertex set into finite sets $\{T_i \st i\in I\}$ such that {\rm (1)} the subgraphs induced by these tiles $T_i$ are connected and all isomorphic to each other; {\rm (2)} the following tiling graph $\widehat \Gamma$ is isomorphic to $\Gamma$: the vertex set is $I$, and $(i,j)$ is an edge of $\widehat\Gamma$ if{f} there is an edge of $\Gamma$ connecting $T_i$ with $T_j$; {\rm (3)} for each $n\geq 1$, there is such a tiling graph $\widehat\Gamma^{n+1}$ on $\widehat\Gamma^{n}$ in such a way that the resulting nested sequence of tiles $T^n(x)\in \widehat\Gamma^n$ containing any fixed vertex $x$ of $\Gamma$ exhausts $\Gamma$. 

Furthermore, $\Gamma$ has a {\bf strongly scale-invariant tiling} if each $T^n$ is isomorphic to $T^{n+1}$. 

If $\Gamma$ has a scale-invariant tiling, is it necessarily of polynomial growth?
\eq

Note that even for a strongly scale-invariant tiling, there could be several ways to iterate the procedure and get $\widehat\Gamma^{n+1}$ from $\widehat\Gamma^{n}$. For instance, on $\Z^d$, the boxes $[0,2^k)^d$ do not exhaust the graph.

As we will explain below, the scale-invariance of a non-Abelian group does not seem to imply alone the existence of a scale-invariant tiling of any of its Cayley graphs. However, from certain expanding homomorphims of the {\bf real Heisenberg group} we will obtain self-similar actions of the integer Heisenberg group with nice properties that imply the existence of a strongly scale-invariant tiling on some Cayley graph of it. (See Theorem~\ref{t.Heis}.) We are unable to generalize this procedure for not virtually nilpotent groups, but, for our amenable examples, we can get a {\bf F{\o}lner monotiling} (i.e., a sequence of tilings $\widehat\Gamma^n=\{ g T_n: g \in G_n\}$, each with a single connected prototile $T_n$ that form a F{\o}lner sequence of $\Gamma$ as $n\to\infty$) with the extra property that $G_n\simeq G$. Typically, these tiling graphs $\widehat\Gamma^n$ seem to grow in some sense, instead of being isomorphic to each other. Still, such a sequence might help in certain weaker versions of renormalization and provide interesting results. 

However, there are some further ingredients that need to be generalized from $\Z^d$. The main probabilistic challenge is the following, related to the conjectural uniqueness of the giant cluster in percolation on finite transitive graphs \cite{ABS}. For definitions, see Section~\ref{s.perc}, where we will discuss the issues related to scale-invariant tilings and percolation.
 
\bq\label{q.unique} Let $\Gamma$ be an amenable transitive graph, and let $\Cl_\infty$ be its unique infinite percolation cluster at some $p>p_c(\Gamma)$, with density $\theta(p)$. For a finite vertex set $W \subset \Gamma$, let $c_i(W)$ denote the number of vertices in the $i^\text{th}$ largest connected component of $W$. Does there exist a connected F{\o}lner sequence $F_n\nearrow \Gamma$ such that for almost all percolation configurations,
$$
\lim_{n\to\infty} \frac{c_2(F_n \cap \Cl_\infty)}{c_1(F_n \cap \Cl_\infty)} = 0\,,
$$
moreover,
$$
\lim_{n\to\infty} \frac{c_1(F_n \cap \Cl_\infty)}{|F_n|}=
\lim_{n\to\infty} \frac{|F_n \cap \Cl_\infty|}{|F_n|}=\theta(p)\,?
$$
\eq

\medskip

Regardless of possible applications to statistical mechanics, it would be interesting to understand the class of scale-invariant and strongly scale-invariant groups better. More generally, given an arbitrary group $G$, one could study its subgroup
$$G_\infty:=\bigcap\Big\{H \; : \; H\leq G,\ [G:H]<\infty,\ G\simeq H\Big\}\,.$$
For what groups is this subgroup trivial or small, say, isomorphic to $\Z$ or Abelian?

Not every nilpotent group is strongly scale-invariant: maybe a bit surprisingly, there exist torsion-free nilpotent groups that are even {\bf co-Hopfian} (i.e., they have no proper subgroups isomorphic to themselves) and hence $G_\infty=G$ \cite{Beleg}. 

There are many other groups that are known to be co-Hopfian \cite[Item III.~22]{dlH}. Z.~Sela proved that a torsion-free non-elementary {\bf hyperbolic group} is co-Hopfian if{f} it is not a non-trivial free product \cite{Sela}. The simplest non-elementary non-co-Hopfian hyperbolic group is the free group $F_r$ (with $r\geq 2$): it has many proper subgroups isomorphic to itself, but none is of finite index,
since $s-1=[F_r:F_s]\,(r-1)$ if $[F_r:F_s]<\infty$. Hence 
$F_r$ is not scale-invariant, either. More generally, if there is an Euler-Poincar\'e characteristic $\chi$ (i.e., $\chi(H)=\chi(G)\,[G:H]$ if $[G:H]<\infty$) that is non-zero for $G$, then $G$ cannot be scale-invariant. An example of such a $\chi$ is the first $\ell^2$-Betti number of the group, which is also the von Neumann dimension of the $G$-invariant Hilbert space of harmonic functions with finite Dirichlet energy, see \cite{Paschke, BV}. Thus, if $G$ has non-constant harmonic Dirichlet functions, i.e., its first $\ell^2$-Betti number is non-zero, then it is not scale-invariant. Since non-trivial free products have infinitely many ends, and thus their first $\ell^2$-Betti number is non-zero \cite{SW, BV}, this discussion establishes the following proposition. It was first suggested to us by Benjamini, but, as we have recently learnt from Sapir, it had been proved earlier in \cite{Bridsonetal}. They also started from Sela's theorem, but instead of using non-constant harmonic Dirichlet functions, they concluded with a topological argument bounding the algebraic rank (the minimal number of generators) of finite index subgroups of non-trivial free products.

\begin{proposition}\label{pr.hyperb}
Torsion-free non-elementary hyperbolic groups are not scale-invariant. \qed
\end{proposition}

More generally, M.~Sapir conjectures \cite{Sapir} that non-elementary relatively hyperbolic groups are not scale-invariant, and suggests that the methods of \cite{DruSap} could work to prove this. What easily follows from \cite[Theorem 1.14]{DruSap} is that any scale-invariant non-elementary relatively hyperbolic group must be a free product amalgamated over a virtually cyclic or parabolic subgroup.

\section{The general construction}\label{s.general}

\noindent{\it Proof of Theorem~\ref{t.general}}.
Consider the right coset tree $\T$ of the subgroup sequence $(H_n)_{n\geq 0}$: the root is $H=H_0$, and a coset $H_{n+1} y$ is a child of $H_n x$ if $H_{n+1} y \subset H_n x$. The number of children of $H_n x$ is $[H_n:H_{n+1}]$. For a ray $H=H_0 x_0 \supset H_1 x_1 \supset H_2 x_2 \supset \dots$ in $\T$ we will use the shorthand notation $x=(x_1,x_2,\ldots)$; the set of these rays is the boundary $\p\T$ of the tree, equipped with the usual metrizable topology. If we have normal subgroups, $H_n\lhd H$ $\forall n$, then $\p\T$ can be equipped with a group structure: it is the profinite completion of $H$ with respect to the series $(H_n)_{n\geq 0}$, see e.g.,~\cite{profinite}.

Since $H_n$ is $A$-invariant, the semidirect product $G:=A\ltimes H$ acts on $\T$ by the affine transformations $(H_n x)^{(\alpha,h)}=H_n\, \alpha(x)h$. Clearly, this action is transitive on the levels of $\T$, and it extends to a continuous action on $\p\T$. The group of affine transformations leaving $H_n$ invariant is the semidirect product $A\ltimes H_n$, and for any given vertex $H_n x$ on the $n$th level of the tree, its stabilizer is
\be\label{e.Stab}
\St_G(H_n x)=\Big\{\big(\alpha,\,\alpha(x)^{-1}h_n x\big) \st \alpha\in A,\ h_n \in H_n\Big\}\simeq A\ltimes H_n\simeq A\ltimes H\,.
\ee
If $H_{n+1} y$ is a child of $H_n x$, then $\St_G(H_{n+1} y)\subset \St_G(H_n x)$ is a subgroup with index $[H_n : H_{n+1}]$. For a point in the boundary, $\overline{x}=(x_1,\ldots) \in \p\T$, we have $\St_G(\overline{x})=\bigcap_{n\geq 1} \St_G(H_n x_n)$. Thus, if we prove that there exists $\overline{x}\in\p\T$ with a finite (or even trivial) stabilizer, then the sequence $G_n:=\St_G(H_n x_n)$ will show that $G$ is scale-invariant.

Fix an element $(\alpha,h)\not=(1,1)$ of $A\ltimes H$. We claim it cannot stabilize every point of an open set in $\p T$. Otherwise, there would be an $x_n$ such that the entire subtree below $H_n x_n$ is stabilized, which, since $\bigcap_{i\geq n} H_i=\{1\}$, means that $H_n x_n$ is pointwise stabilized by $(\alpha,h)$. That is, $\alpha(h_n x_n)h=h_n x_n$ for all $h_n\in H_n$. In particular, $\alpha(x_n)h=x_n$, hence $\alpha(h_n)=h_n$, which is possible only if $\alpha=1$ because the action of $A$ on $H_n$ is supposed to be faithful. Then we immediately get $h=1$, too.

On the other hand, since $(\alpha,h)$ acts continuously on $\p\T$, the set $S_{(\alpha,h)}\subset \p\T$ of points stabilized by it is closed. Together with the previous paragraph, the complement $S_{(\alpha,h)}^c$ is open and dense in $\p\T$. Since $\p\T$ is a compact metrizable space and $A\ltimes H$ is countable, by Baire's category theorem we have that $\bigcap \Big\{ S_{(\alpha,h)}^c \st (\alpha,h)\in (A\ltimes H) \setminus \{(1,1)\} \Big\}$ is nonempty. In other words, there is some $\overline{x} \in\p\T$ whose stabilizer is the trivial $\{(1,1)\}$, and we are done.\qed
\bigskip

We will now discuss two issues that turn out to be related to each other. Firstly, we would like to see that the conditions of our theorem are fulfilled sometimes. Secondly, we would like to describe our above isomorphisms $G\longrightarrow G_n$ a bit more explicitly, e.g., to see if they could have a subsequence $(n_k)$ along which they are equal to $\phi^{\circ k}$ for some $\phi$, thus proving strong scale-invariance.

The next lemma describes a natural situation in which the conditions on how $A$ should act on each $H_n$ hold automatically. The straightforward proof is left to the Reader.

\bl\label{l.twist} Consider a family of injective endomorphisms $\{\psi_i\st i\in I\}$ of $H$ such that for each $i\in I$, for some $\tau_i\in \Aut(A)$ and all $h\in H$ and $\alpha\in A$, we have
\be
\psi_i(\alpha(h)) = \tau_i(\alpha)(\psi_i(h))\,.\label{e.twist}
\ee
Let $H_n=\psi_{i_1}\circ\ldots\circ\psi_{i_n}(H)$ for some $i_1,\ldots,i_n \in I$. Then $A$ acts faithfully on $H_n$, and $A\ltimes H_n \simeq A\ltimes H$.
\el

In all our examples, the group $H$ will be Abelian, and the set $I$ in the previous lemma will be a singleton, i.e., we will have $H_n=\psi^{\circ n}(H)$ for all $n\geq 0$. Then, as we will see in Proposition~\ref{pr.homos} in a second, for the sequence $(G_n)_{n\geq 0}$ we have constructed above, the isomorphisms $G \longrightarrow G_n$ do arise from composing a finite set of endomorphisms $J$ in the way described in the lemma, with $J$ as $I$. This gives hope that our construction could be iterated to obtain more scale-invariant groups: we would just need a non-trivial group $A$ of automorphisms of $G$ that satisfy~(\ref{e.twist}) with this finite set $J$. Unfortunately, we do not know if such an $A$ exists for any of our examples below.

\bpr\label{pr.homos} 
{\bf (i)} In the setting of Theorem~\ref{t.general}, assume furthermore that $H_n=\psi^{\circ n}(H)$ for some injective endomorphism $\psi$. Then, for any vertex $v$ in the $(H_n)$-coset tree $\T$, the stabilizer $\St_G(v)$ is of the form $\phi_{v}(G)$, and if the parent of $v$ is $u$, then $\phi_v=\phi_u\circ \phi_{u,v}$, where $\phi_{u,v} \in J$, for a finite set $J$ of injective endomorphisms of $G$ with $|J|=[H:\psi(H)]$.\\
{\bf (ii)} In particular, if there is a periodic ray $\overline{x}=(x_1,x_2,\dots) \in\p\T$ (i.e., $\big(\phi_{v_n,v_{n+1}}\big)_{n=0}^\infty$ is periodic on $J$, where $v_n=(x_1,\dots,x_n) \in \T$) such that $\St_G(\overline{x})$ is finite, then $G$ is strongly scale-invariant.
\epr

\proof
Let  $y_1,\dots,y_t$ be a set of right coset representatives for $H_1=\psi(H)$ in $H$. Then, by (\ref{e.Stab}), we have the isomorphism $\phi_i:\, G\longrightarrow \St_G(H_1 y_i)$, for each $i=1,\ldots,t$, given by 
\be\label{e.phii}
\phi_i(\alpha, h):= (\alpha,\,\alpha(y_i)^{-1}\psi(h)y_i)\,.
\ee
We will have $J=\{\phi_1,\dots,\phi_t\}$. Indeed, for $H_2=\psi^{\circ 2}(H)$, a set of right coset representatives in $H$ is $\psi(y_j)y_i$, with $i,j=1,\ldots,t$, and we have the isomorphisms $\phi_{ij}:\, G\longrightarrow \St_G(H_2\, \psi(y_j)y_i)$ given by
\begin{align*}
\phi_{ij}(\alpha, h):&= \Big(\alpha,\,\alpha\big(\psi(y_j)y_i\big)^{-1}\,\psi^{\circ 2}(h) \,\psi(y_j)y_i\Big)\\
&= \phi_i \big(\alpha,\,\alpha(y_j)^{-1}\psi(h)y_j\big)=\phi_i\circ \phi_j \, (\alpha,h)\,.
\end{align*}
Continuing by induction, we get that $H_n=\psi^{\circ n}(H)$ has a set of right coset representatives
$$
x_{i_1\cdots i_n}:=\psi^{\circ (n-1)}(y_{i_n})\cdots \psi(y_{i_2})\, y_{i_1}
$$
with $i_k=1,\ldots, t$ for each $k=1,\ldots, n$. Note here that the ray leading from the root to $H_n x_{i_1\cdots i_n}$ is $(x_{i_1},\ldots, x_{i_1\cdots i_n})$. Denoting $G(i_1,\ldots,i_n):=\St_G(H_n x_{i_1\cdots i_n})$ and $G(\emptyset)=\St_G(H)=G$, the isomorphisms $\phi_{i_1\cdots i_n}:\, G\longrightarrow G(i_1,\ldots,i_n)$ then satisfy
\be\label{e.comp}
\phi_{i_1\cdots i_n}=\phi_{i_1}\circ\ldots\circ \phi_{i_n}\,.
\ee
Note the order of composition: the isomorphism $G(i_1,\ldots,i_{n-1})
\longrightarrow G(i_1,\ldots,i_n)$ is {\it not} $\phi_ {i_n}$.

What would ensure that $\bigcap_{k\geq 0} G(i_1,\ldots,i_k)=\bigcap_{j\geq 0}\phi^{\circ j}(G)$ for a single injective endomorphism $\phi$\,? The only reasonable answer seems to be that the infinite sequence $\overline{i}:=(i_1,i_2,\ldots)$ should be periodic: if there exists some $p\geq 1$ with $i_{k+p}=i_k$ for all $k\geq 1$, then $\phi:=\phi_{i_1\cdots i_p}$ would do. \qed
\medskip

So we have arrived at the question: {\it is there an infinite periodic ray $\overline{x}\in\T$ whose $G$-stabilizer is trivial, or finite, at least?} As we will see, the answer is negative in all the scale-invariant cases we have analyzed. Of course, this does not prove that these examples are not strongly scale-invariant, but we have no further ideas to attack this problem.

\section{The examples}\label{s.examples}

\subsection{The lamplighter groups}\label{ss.LL}

For simplicity, we first discuss the case $\F=\Z_2$.

Let $H$ be the additive subgroup of the group $\Z_2[[t]]$ of formal power series over $\Z_2$ consisting of finite Laurent polynomials of $(1+t)$, and consider the injective endomorphism $\psi(F(t)):=tF(t)$ for $F(t)\in \Z_2[[t]]$. Since $tF(t)=(1+t)F(t)-F(t)$, we have that $\psi(H)\subseteq H$. Observe that $(1+t)^k-1 \in \psi(H)$ for any $k\in\Z$; this easily implies that $\psi(H)$ is exactly the subgroup of $H$ of power series divisible by $t$, with index $[H:H_1]=2$. We then let $H_n:=\psi^{\circ n}(H)$, a nested sequence of finite index isomorphic subgroups. The boundary $\p\T$ of the coset tree is the
profinite additive group $\Z_2[[t]]$, via the identification
\be\label {e.Phi}
\Phi:\, x_1x_2\dots \mapsto \sum_{i\geq 1} x_i t^{i-1}\,.
\ee

Now let $A$ be the cyclic group $\Z$ acting on $H$ by multiplication by $(1+t)$. Thus the semidirect product $\G=A\ltimes H$ is the group of the following transformations of $\Z_2[[t]]$:
\be\label{e.trafo}
F(t) \mapsto (1+t)^m F(t) + \sum_{k\in\Z} f(k) (1+t)^k\,,
\ee
where $m\in \Z$ and $f:\,\Z\longrightarrow \Z_2$ is any function with finitely many non-zero values.

This group $\G=\big(\oplus_{\Z}\, \Z_2\big) \rtimes \Z = \Z_2 \wr \Z$ is the standard lamplighter group; for each element $(m,f)$, one can think of $m\in \Z$ as the position of the lamplighter, while $f:\,\Z\longrightarrow \Z_2$ is the configuration of the lamps. We will sometimes represent $f$ by the finite set $\supp f\subset \Z$. The usual wreath product generators are $s$ and $R$, representing ``switch'' and ``Right''; we will also use $L=R^{-1}$. So, for example, $Rs=(1,\{1\})$. In terms of the representation~(\ref{e.trafo}), the action of $s$ is $F(t)\mapsto F(t)+1$, while the action of $R$ is $F(t)\mapsto (1+t)F(t)$.

Since our $\psi:H\longrightarrow H$ clearly commutes with the action of $A$, we can apply Lemma~\ref{l.twist}, and Theorem~\ref{t.general} shows that the lamplighter group $\G=A\ltimes H$ is scale-invariant.
\medskip

The action of the lamplighter group on the infinite binary tree $\T$ can now be described by the combination of~(\ref{e.Phi}) and~(\ref{e.trafo}), and it turns out to be a finite-state self-similar action. We recall now the basic definitions, but see \cite{GNS,fractal,Nekra} for further details and background. 

\bde\label{d.ss}
The action of a group $G$ on the $b$-ary tree $\T_b$ ($b\geq 1$) is called {\bf self-similar} if for any $g\in G$, any letter $x\in \{0,1,\dots,b-1\}$, and any finite or infinite word $w$ on this alphabet, there is a letter $y$ and $h\in G$ such that $(xw)^g=y(w^h)$. If $S\subseteq G$ generates $G$ as a semigroup, and $\forall s\in S$ and word $xw$ there is a letter $y$ and $t\in S$ such that $(xw)^s=y(w^t)$, then $S$ is called a {\bf self-similar generating set}. Then the group can clearly be generated by an automaton with states $S$. The usual diagram of this automaton is called the {\bf Moore diagram} of $S$. If there is a finite such $S$, then $G$ is called a {\bf finite-state self-similar group}.  
\ede

For a self-similar action by $G$, for any $g\in G$ and finite word $v$ there is a word $u$ of the same length and $h\in G$ such that $(vw)^g=u(w^h)$ for any word $w$. This  $h$ is called the {\bf restriction} $h=g|_v$, and we get an action of $G$ on the subtree starting at $v$. The action of the full automorphism group of $\T_b$ is of course self-similar, and there is  the obvious wreath product decomposition 
\be\label{e.treathA}
\Aut(\T_b) \simeq \Aut(\T_b) \wr \mathrm{Sym}_b\,,
\ee
corresponding to the restriction actions inside the $b$ subtrees at the root and then permuting them. For a general self-similar action by $G\leq \Aut(\T_b)$, the isomorphism~(\ref{e.treathA}) gives an embedding
\be\label{e.treathG}
G  \hookrightarrow  G \wr \mathrm{Sym}_b\,.
\ee

Coming back to the lamplighter group $\G$, its self-similarity was first noticed and proved by Grigorchuk and \.Zuk in \cite{GZ}, but the above representation using $\Z_2[[t]]$ gives a much simpler proof, found by \cite{GNS}. Namely, consider the following new generators of the lamplighter group: $a=Rs$, $b:=R$. Note that $s=b^{-1}a=a^{-1}b$. Then, the action of these generators on the binary tree $\T$ can be easily checked to be
\be\label{e.GZ}
(0w)^a = 1w^b \qquad & & \qquad (0w)^b = 0w^b\\
(1w)^a = 0w^a \qquad & & \qquad (1w)^b = 1w^a,\nonumber
\ee
for any finite or infinite $\{0,1\}$ word $w$. Hence $\{a,b\}$ is a finite self-similar generating set.
Another usual notation for this self-similar action, using~(\ref{e.treathG}), is
\be\label{e.selfsim}
\qquad a=(b,a)\eps\,, \hskip 2 cm b=(b,a)\,,
\ee
where $(g,h)$ is the tree-automorphism acting like $g$ on the 0-subtree and like $h$ on the 1-subtree, $\eps$ is switching these two subtrees, and the order of the multiplication is dictated by having a right action on the tree. In particular, $(g,h)(g',h')=(gg',hh')$ and $(g,h)\eps=\eps(h,g)$.

We note that in the literature there are a few slightly different versions of (\ref{e.selfsim}) to describe the lamplighter group. This is partly due to the fact that interchanging the generators $a$ and $b$ induces an automorphism $\iota$ of $\G$, see e.g.~\cite{GZ}.

Let us now see what the endomorphisms~(\ref{e.phii}) are in terms of the automaton representation.

\bpr\label{pr.homosLL}
A suitable set $J$ of endomorphisms in part (i) of Proposition~\ref{pr.homos} for the lamplighter group, with the notation of~(\ref{e.selfsim}) and the automorphism $\iota$, is given by
\be\label{e.pgi}
\phi_0(g)=\big(g , \; \iota (g) \big) \qquad\textrm{and}\qquad
\phi_1(g)=\big(\iota (g), \; g \big).
\ee
\epr

\proof With the coset representatives $y_0=0$ and $y_1=1$ for $H_1$ in $H$, writing $(m,f)\in \G$ as $(m,h(t))$ via $h(t)=\sum_{k\in\Z} f(k)(1+t)^k$, we get
$$
\phi_0\big(m,h(t)\big) = \big(m,\,th(t)\big)\qquad\textrm{and}\qquad
\phi_1\big(m,h(t)\big)= \big(m,\,th(t)+(1+t)^m+1\big).
$$
In particular, for $a=(1,1)$ and $b=(1,0)$, we get $\phi_0(a)=\phi_1(b)=(1,t)=a^{-1}ba=sRs$ and $\phi_1(a)=\phi_0(b)=b$. 

To prove now (\ref{e.pgi}), it is enough to check it for the generators $a,b$. For $\phi_1(a)=\phi_0(b)=b$, this is trivial from~(\ref{e.selfsim}). For $\phi_0(a)=\phi_1(b)=a^{-1}ba$, we have
$$
a^{-1}ba=\eps (b^{-1},a^{-1})(b,a)(b,a)\eps=(a^{-1}aa,\; b^{-1}bb)=\big(a , \; \iota (a) \big)=\big(\iota(b), \; b \big),
$$
and~(\ref{e.pgi}) is proved.\qed 
\medskip

The form (\ref{e.pgi}) shows it clearly that the $\phi_i$ are isomorphisms, and inductively, that for any finite word $\overline{y}=x_1\cdots x_n\in\T$, the stabilizer $\St_\G(\overline{x})$ equals $\phi_{x_1}\circ\dots\circ\phi_{x_n}(\G)$, as we also showed in~(\ref{e.comp}).
\medskip

We now discuss which rays $\overline{x}\in\p\T$ can have a trivial stabilizer, and show in particular that part (ii) of Proposition~\ref{pr.homos} cannot imply that $\G$ is strongly scale-invariant.

Recall the map $\Phi$ from (\ref{e.Phi}). If, for some $\overline{x}\in\p\T$, there is a non-trivial element $g\in \St_\G(\overline{x})$, then we have a finite sequence of transformations $F(t)\mapsto F(t)\pm 1$ and $F(t)\mapsto (1+t)^{\pm 1} F(t)$ fixing the power series $\Phi(\overline{x})\in \Z_2[[t]]$. This implies that $\Phi(\overline{x})$ is a rational function $U(t)/V(t)$ with $U(t),V(t) \in \Z_2[t]$, where
$U(t)=(1+t)^{\ell_n}+\dots+(1+t)^{\ell_1}$ and $V(t)=(1+t)^{\ell_0}+1$, with $\ell_i\in\Z$.
There are only countably many such functions, while continuum many possible words $\overline{x}$, so for most words $\overline{x}$ we have $\St_\G(\overline{x})=\{1\}$. But can we have $\St_\G(\overline{x})=\{1\}$ for a periodic word $\overline{x}:=\overline{y}\,\overline{y}\cdots$, with $\overline{y}=x_1\cdots x_k$? Note that the corresponding power series is $\Phi(\overline{x})=(\sum_{i=1}^{k} x_i t^{i-1})(1+t^k+t^{2k}+\dots)$.

\bpr\label{pr.periodLL} 
For any periodic ray $\overline{x}$, the stabilizer $\St_\G(\overline{x})$ is infinite.
\epr

\proof
We first show that there is a non-trivial element of $\G$ that fixes $\Phi(\overline{x})$.
Using our above representation of fixed points by $U(t)/V(t)$, and noticing that the finite combinations of functions $(1+t)^{\ell_i}$ with $\ell_i\in\N$ are exactly all the polynomials in $\Z_2[t]$, one can immediately translate our claim into the following:

\bl\label{l.perioLL}
Given any $k\in\Z^+$, there exist integers $0\leq \ell < m $ such that $\big( (1+t)^{m}+(1+t)^{\ell}\big)\,\big(1+t^k+t^{2k}+\dots\big)$ is a polynomial in $\Z_2[t]$.
\el

\proof The coefficient of $t^j$ in $\big( (1+t)^{m}+(1+t)^{\ell}\big)\,\big(1+t^k+t^{2k}+\dots\big)$ is
$$
{m \choose j} + {\ell \choose j} + {m \choose j-k} + {\ell \choose j-k} +
{m \choose j-2k} + {\ell \choose j-2k} + \dots .
$$
We want this to be zero (mod 2) for all large enough $j$. This is equivalent to
$$
F_k(m,i)+F_k(\ell,i)=0 \quad (\mbox{mod }2),\qquad\mbox{where}\qquad
F_k(m,i):=\sum_{r=0}^{\lfloor\frac{m-i}{k}\rfloor} {m \choose i+rk}\,,
$$
for all $i=0,\dots,k-1$. Lucas' theorem (1878) on the parity of binomial coefficients says that $a\choose b$ is odd if{f} each binary digit $a_i$ of $a$ is larger than the corresponding digit $b_i$ of $b$, see e.g.~\cite{binomial}. In particular, if $m=2^t-1$ for some $t\in\Z^+$, then each term in $F_k(m,i)$ is odd, hence $F_k(m,i)=\lfloor\frac{m-i}{k}\rfloor + 1$ (mod 2). Therefore, if we find $m=2^t-1$ and $\ell=2^s-1$ such that $m=Mk+a$ and $\ell=Lk+a$ with $M,L,a\in \N$, then $F_k(m,i)+F_k(\ell,i)=M+L$ (mod 2) for each $i$.

Let us write $k=2^\kappa K$ with $\kappa\in\N$ and $K$ odd, and take $a=2^\kappa$-1; then we want two different integers $M=(2^t-2^\kappa)/(2^\kappa K)=(2^{t-\kappa}-1)/K$ and $L=(2^{s-\kappa}-1)/K$, such that $M+L$ is even. If $M,L$ given by the above formulas are integers, then they are necessarily odd, so we just need that $K \, | \, 2^u-1$ for two different integers $u$. Since $K$ is odd, by the Euler-Fermat theorem we have $K \, | \, 2^{v\phi(K)}-1$ for any $v\in\Z^+$, and we are done.\qed
\medskip

Thus, any periodic ray has a non-trivial stabilizer. This easily implies that the stabilizer is in fact infinite.
If $(\alpha,h)\in A\ltimes H$ stabilizes $\overline{x}\in\p\T$, then $(\alpha,h)^k$ also does, for any $k\in \Z$. If the automorphism group $A$ is torsion free (as is the case now), then all elements $(\alpha,h)^k$ are different, provided $\alpha\not=1$; hence $\St_\G(\overline{x})$ is infinite. If $\alpha=1$, and $H$ is Abelian (again, as is the case now), then (\ref{e.Stab}) implies that the stabilizer $\St_\G(\overline{x})$ contains the entire subgroup $(1,H)$, and we are done.\qed
\medskip

\noindent{\bf Remark.} Actually, our subgroup chains $(G_n)_{n\geq 0}$ already appear in \cite{GZ}, where the self-similar action was used for studying simple random walk on $\G$. For this, they needed a bounded index subgroup chain with trivial intersection, and this chain was realized as a nested sequence of vertex stabilizers, using a Baire category argument similar to ours. Moreover, they proved that the stabilizer of each boundary point is either trivial or isomorphic to $\Z$. However, the inner structure of these stabilizers was not important for them (besides being able to continue the subgroup chain). In particular, they did not observe the existence of the isomorphisms~(\ref{e.pgi}).

Accordingly, the key observation that started our work was ``orthogonal'' to~\cite{GZ}: namely, that the {\bf Diestel-Leader graph} $DL(2,2)$ is the Cayley graph of $\G$ with the generators $\gen{R,Rs}$ on one hand, and of the index two subgroup $G_1=\gen{Rs,sR}$ on the other. Our results above show that this isomorphism of the Cayley graphs is due to the group isomorphism $g\mapsto \iota(\phi_0(g))$. See \cite{DL,W:DL,ariel,EFW} for the definition of these graphs and some background.
\medskip

For the general case, when $\F$ is an arbitrary finite Abelian group, one can use its decomposition as a direct sum of cyclic groups, and then $\F$ and $\F [[t]]$ can again be considered as rings. This was used by \cite{SiSt} to describe $\F\wr\Z$ as a group of transformations of $\F[[t]]$, which showed that $\F \wr \Z$ is again a group of a finite automaton. Then, everything we did above goes through in this generality.

\subsection{The solvable Baumslag-Solitar groups}\label{ss.BS}

Let $m,\ell$ be a pair of positive integers that are relative primes. Let $H$ be the additive group of rational numbers of the form $a/m^b$ with $a,b\in \Z$, and $\psi:\,H\longrightarrow H$ be multiplication by $\ell$, an injective endomorphism. Then $H_n:=\psi^{\circ n}(H)$ is the subgroup of rationals of the form $\ell^n a/m^b$. The associated coset tree $\T$ is $\ell$-ary, and $\p\T$ can be identified with the profinite group of $\ell$-adic integers $\Z_{(\ell)}$. Now, the cyclic group $\Z$ acts on $H$ by multiplication by $m^t$, $t \in \Z$, and this action commutes with $\psi$, hence Lemma~\ref{l.twist} and Theorem~\ref{t.general} apply, and $G=\Z \ltimes H$ is a scale-invariant group.

The transformations $s_i:\, u\mapsto mu+i$ on $\Z_{(\ell)}$ for $i=0,1,\dots,m-1$ generate the action of $G$ on $\T$, and this action is self-similar:
$$
(jw)^{s_i} = \big(mj+i \mbox{ (mod }\ell)\big)\, w^{s_{\lfloor(mj+i)/\ell\rfloor}},
$$
where $j$ is a letter and $w$ is a finite or infinite word in $\{0,1,\ldots,\ell-1\}$.
Taking $b:=s_0$ and $a:=s_1s_0^{-1}$, one can show that $G$ is in fact the
Baumslag-Solitar group $\BS(1,m)=\Gen{a,b \mid bab^{-1}=a^m}$. This representation of $\BS(1,m)$ was discovered by Bartholdi and \v{S}uni\'k \cite{solvauto}.

Let us show that the stabilizer of a periodic ray in $\T$ could never be finite; hence Proposition~\ref{pr.homos} cannot show that $\BS(1,m)$ is strongly scale-invariant. Rewriting the vertex stabilizers (\ref{e.Stab}) for the present case, when the sequence $x_n$ of coset representatives is periodic with some period $p$, we want to show that for any $p\in\Z^+$, $a\in \{0,1,\dots,\ell^p-1\}$, $b\in\Z$, there exists some non-zero $t\in\Z$ such that
$$
\bigcap_{k=1}^\infty \left\{
(1-m^t)\frac{a}{m^b}\big(1+\ell^p+\ell^{2p}+\dots+\ell^{(k-1)p}\big)+\frac{\ell^{kp}\Z}{m^{\Z}}
\right\} \not= \emptyset\,.
$$
Note that $1+\ell^p+\ell^{2p}+\dots+\ell^{(k-1)p}=(\ell^{kp}-1)/(\ell^p-1)$, and write $\ell^p-1=m^c d$ with $c,d\in\N$, where $(m,d)=1$. By the Euler-Fermat theorem, we can choose $t\in\Z^+$ such that $m^t-1=df$ with $f\in \Z^+$. This way our intersection becomes
$$\bigcap_{k=1}^\infty \left\{
\frac{-fa}{m^c}\big(\ell^{kp}-1\big)+\frac{\ell^{kp}\Z}{m^{\Z}}
\right\} \ni \frac{-fa}{m^c} \,,
$$
and we are done. From the stabilizer being non-trivial, we immediately get that it is in fact infinite, by the argument at the end of Proposition~\ref{pr.periodLL} of the lamplighter case.

\subsection{The affine groups}\label{ss.affine}

Let $H$ be the additive group $\Z^d$, and $H_n$ be the subgroup $2^n\Z^d$. The action of $GL(\Z,d)$ commutes with the multiplication endomorphisms $Z^d \longrightarrow 2^n\Z^d$. Hence, by Lemma~\ref{l.twist} and Theorem~\ref{t.general}, $G=GL(\Z,d)\ltimes \Z^d$ is a scale-invariant group. This action of $G$ on the $2^d$-ary tree $\T$ is in fact self-similar, and it was first discovered by Brunner and Sidki \cite{affine}.

Again, the stabilizer of a periodic ray in $\T$ is never trivial. This time, the reason is that we need that for any $p\in\Z^+$ and $v\in \{0,1,\dots,2^p-1\}^d$ there exists an $\alpha\in GL(\Z,d)\setminus\{I\}$ such that
\be\label{e.intersect}
\bigcap_{k=1}^\infty \left\{
(I-\alpha) v \big(1+2^p+2^{2p}+\dots+2^{(k-1)p}\big)+2^{kp}\Z^d
\right\} \not= \emptyset \,,
\ee
and this follows immediately from the next lemma.

\bl\label{l.perioGL} Let $d\geq 2$. For any $v\in \Z^d$ there exists $\alpha\in GL(\Z,d)\setminus\{I\}$ with $\alpha v=v$.
\el

\proof For $d=2$, given $v=(x,y)^T\not=(0,0)^T$, it is easy to find the following solution $\alpha=\alpha(x,y)\in GL(\Z,2)\setminus\{I\}$:
$$
\left(\begin{array}{cc}
xy+1 & -x^2 \\
y^2 & -xy+1\\
\end{array}\right)
\left(\begin{array}{c}
x \\
y \\
\end{array}\right)
=
\left(\begin{array}{c}
x \\
y \\
\end{array}\right).
$$
For $v=(0,0)^T$ any $\alpha$ will do. Now, for $v=(x_1,\ldots,x_d)^T \not= \overline{0} \in \Z^d$, with $d\geq 3$, we can assume by permuting the coordinates that $x_1\not=0$, and then just use the matrix with $\alpha(x_1,x_2)$ in its upper right corner and 1's along the diagonal from the third entry downward.
\qed
\medskip

Again, we easily get that the stabilizer is infinite. For each $k\in \Z$, we have $\alpha^k v=v$ for the above $\alpha$, and all the matrices $\alpha^k$ are different, since $\alpha^k=k(\alpha-I)+I$, as it is easy to check.
\medskip

For $G=A \ltimes \Z^d$ with any subgroup $A \leq GL(\Z,d)$, the above use of Lemma~\ref{l.twist} and Theorem~\ref{t.general} goes through, hence these $G$ are again scale-invariant. However, Lemma~\ref{l.perioGL} does not apply in general to show that we do not get strongly scale-invariant subgroup chains for them, since there are many subgroups $A$ without non-trivial elements fixing a given integer vector. Still, we expect that periodic rays will always have infinite stabilizers. We prove this only for the case when $A$ contains a non-trivial unipotent matrix; in particular, whenever $G=A \ltimes \Z^d$ is virtually nilpotent. For semidirect products given by hyperbolic matrices, such as the group $Sol$, one reason to believe that the stabilizers will be infinite is that their geometry is quite similar to that of the Diestel-Leader graphs, as demonstrated in \cite{EFW}.

The next lemma clearly implies our claim above on subgroups $A$ containing unipotent matrices:

\bl\label{l.unipot} 
For any integer-valued unipotent matrix $M\in GL(\Z,d)$, any $p\in\Z^+$ and $v\in \{0,1,\dots,2^p-1\}^d$, there exists an $n\in\Z^+$ such that the intersection~(\ref{e.intersect}) is non-trivial for $\alpha=M^n$. 
\el

\proof
It is easy to show (e.g., by induction on the size of the matrix) that $M$ can be written as $I+S^{-1}NS$, where $S\in GL(\Z,d)$ and $N$ is a strictly upper-triangular integer-valued matrix. Hence $M^n=S^{-1}(I+N)^n S$. Note that $(I+N)^n$ has entries above the diagonal which are all combinations of monomials in the entries of $N$
with binomial coefficients ${n\choose 1}, \ldots, {n\choose d}$. So, if $n$ is
divisible by $d!(2^p-1)$, then all the off-diagonal entries of $(I+N)^n$ will be divisible by
$2^p-1$, hence all the entries of $M^n-I=S^{-1}((I+N)^n-I)S$ will be divisible by $2^p-1$. Therefore, with $\alpha=M^n$, the intersection (\ref{e.intersect}) becomes $\bigcap_{k\geq 1} \big\{ \widetilde{M} v (2^{pk}-1)+2^{pk}\Z^d \big\}$ with some integer matrix $\widetilde M$. This intersections contains $-v$, and we are done.
\qed

\section{Percolation renormalization and scale-invariant tilings}\label{s.perc}

Let us start with a rough sketch of how {\bf percolation renormalization} on $\Z^d$ typically works. For further background on percolation, see \cite{Grimm, LPbook}. Consider Bernoulli site percolation on $\Z^d$ with density $p\in [0,1]$, i.e., delete each vertex of $\Z^d$ with probability $1-p$, independently from each other. The connected components of the remaining random graph are called percolation clusters. The probability that the cluster of a given vertex is infinite is denoted by $\theta(p)$, and 
$p_c(\Z^d):=\inf\{p:\theta(p)>0\}$. A basic result, true for percolation on any transitive amenable graph, is that for $p>p_c$, there is almost surely a unique infinite cluster, and its density, measured along any F{\o}lner exhaustion, equals $\theta(p)$. However, the following is much harder to prove, and has been established only for $\Z^d$: 

\bth[Renormalization Lemma, \cite{chemical}]\label{l.AP} On $\Z^d$, $d\geq 2$, for any $p>p_c$ and $\eps>0$, if we take a large box with side-length $n>n_0(p,\eps)$, then with probability at least $1-\eps$ the box is ``good'':  it has a cluster connecting all $2d$ faces of the box (called the giant cluster), while all other clusters in the box have diameter at most $\eps n$. 
\eth

Now, take the lattice tiling by side-length $n$ boxes, then magnify each box from its center by a factor of $5/4$, so that they will slightly overlap. Note that if there are two overlapping $n$-boxes that are both good in a given percolation configuration, then their giant clusters must in fact be connected inside the union of the two boxes (we only need that $\eps<5/8$). Therefore, a large cluster of good $n$-boxes represents well the structure of a large original cluster on scales larger than $n$. If we take $n>n_0(p,\eps)$, the density of good boxes is at least $1-\eps$, and the goodness of two boxes are independent once they do not overlap. Most things about supercritical percolation are easy to understand at density close enough to 1 (using the so-called {\bf Peierls method}, a simple counting argument), even in this slightly dependent case, and by the above argument this knowledge can be transferred to any $p>p_c$. See \cite{Grimm} for more details.

This type of renormalization has been mainly used in showing that the supercritical phase on $\Z^d$ is well-behaved in several different senses, such as: possibly different critical points actually coincide \cite{GrM}, and the large-scale geometry (e.g., length of geodesics, isoperimetric and random walk properties) of the unique infinite cluster is very close to that of $\Z^d$ \cite{chemical, repulsion}. Ideally, one would also like to gain information about behavior at criticality; first of all, to show that critical percolation almost surely has no infinite clusters. However, this is unknown even on $\Z^d$ with $3\leq d \leq 18$; the renormalization method has been enough to show this claim only for percolation in the half-space, and also to show that the half-space has the same critical probability as the entire graph \cite{half}. Let us point out that the failure of this method to understand criticality might be because it does not use the exact scale-invariant structure: it is not important in the method that the graph of the slightly overlapping large boxes is not exactly $\Z^d$ again, neither that it is a somewhat similar graph, only that we can understand very supercritical percolation on it.

A different approach that is related to the scale-invariance of $\Z^d$ is the {\bf renormalization group method}, which has been very useful for several statistical physics models at criticality; see \cite{Brydges}. We will not discuss this approach here.
\medskip

For a general group $G$, what kind of {\bf nice sequence of tilings} would we like to produce in some Cayley graph $\Gamma$ of $G$? Ideally, the following three properties should hold:

\begin{enumerate}
\hitem{i.1}{(1)} The tiling sequence should be scale-invariant in the sense of Question~\ref{q.tiling}.
\hitem{i.2}{(2)} The tiles $T_n$ themselves should be large pieces of $\Gamma$ that ``represent well'' the infinite graph. A natural definition is the Benjamini-Schramm {\bf random local} (or {\bf random weak}) {\bf convergence} \cite{BS:limit,urn}: the finite graphs $\Gamma_n$ converges to a transitive graph $\Gamma$ if for any large $R>0$ and small $\eps>0$, if $n>n_1(R,\eps)$, then at least $1-\eps$ proportion of the vertices of $\Gamma_n$ have their $R$-neighborhood isomorphic to the $R$-ball of $\Gamma$. When the $\Gamma_n$ are transitive, this is equivalent to taking $\eps=0$, and is usually called {\bf local convergence}. If $\Gamma$ is amenable, then any F{\o}lner sequence will converge to it in the random local sense.
\hitem{i.3}{(3)} For amenable $\Gamma$, the tiles should be F{\o}lner sets that have a chance to be intersected by the unique infinite cluster ``substantially'', or more precisely, that have a unique giant cluster with large probability, similarly to the $\Z^d$ case, Theorem~\ref{l.AP}. We do not know exactly what the best general definition of ``goodness'' would be, but cluster sizes are certainly of key importance. So, the tiles should form a F{\o}lner sequence as asked for by Question~\ref{q.unique}.
\end{enumerate}
It is easy to see that for a non-amenable Cayley graph $\Gamma$, no finite {\it subgraphs} $\Gamma_n$ can converge to it in the random local sense. (If there is {\it some} sequence of finite approximating graphs, the group is called sofic. It is not known if there are non-sofic groups \cite{Pestov,urn}.) It is also unclear what the analogue of Property~\iref{i.3} should be. So, we do not presently see how percolation renormalization should work in the non-amenable case.

We now collect what we know which groups have tilings satisfying any of these three properties.
\medskip

First of all, how much the scale-invariance of a group $G$ helps with Property \iref{i.1}, i.e., with producing a sequence of scale-invariant tilings? Unfortunately, in the case of $\Z^d$, the nice isomorphic tiling by large boxes is due not only to the scale invariance, but also to the commutativity of the group. In a general scale-invariant group $G$, with an isomorphic subgroup $H$ of finite index $t$, and a {\it right} Cayley graph $\Gamma(G,S)$ given by a finite symmetric generating set $S$ (i.e., $g$ is a neighbor of $gs$ for any $s\in S$), we have two options to start with.

Probably the better choice is to take a set of {\it right} coset representatives $C=\{g_1,\dots,g_t\}$, i.e., the cosets are $Hg_i$. Then $\{hC : h\in H\}$ is a disjoint covering of $G$, such that the tiles $hC$ are all isometric to each other in the graph metric of $\Gamma$. Let us presently assume that these tiles are connected subgraphs of $\Gamma$. (Or maybe even nice in the sense given above.) Now take the graph $\widehat\Gamma$ induced by this tiling, as defined in Question~\ref{q.tiling}. Unfortunately, there is no reason why the natural left action by $H\simeq G$ on the tiles should imply that this graph $\widehat \Gamma$ is isomorphic to $\Gamma$. (A simple non-Abelian example where the tiling graph is nevertheless isomorphic to the natural Cayley graph of the subgroup is $[F(a,b): F(a^2,b,b^{-1}ab)]=2$, free groups on 2 and 3 generators.)

The other option would be to start with a set $C$ of {\it left} coset representatives. Then $\{Ch : h\in H\}$ is a disjoint covering of $\Gamma$, and the tiling graph $\widehat\Gamma$ is actually a Cayley graph of $H\simeq G$, isomorphic to $\Gamma$. However, the tiles now are not at all isometric to each other inside $\Gamma$, so even if one tile induces a connected subgraph, the vertices in other tiles could be very far from each other. Hence this tiling does not seem to be useful for percolation renormalization in any way.
\medskip

Nevertheless, there exist non-Abelian Cayley graphs with scale-invariant tilings: Theorem~\ref{t.Heis} below will give two such Cayley graphs of the discrete 3-dimensional Heisenberg group. However, these will be Cayley graphs with rather special generating sets, and the construction will use not only the strong scale-invariance of the Heisenberg group, but also that this strong scale-invariance comes in fact from expanding endomorphisms.
\medskip

In general, when the right coset tiling graph $\widehat\Gamma$ is not isomorphic to $\Gamma$, it could still be useful in percolation. It is another Cayley graph of $G$, usually ``larger'' than $\Gamma$, hence the value of the percolation parameter $p$ that is already close enough to 1 to use the Peierls method should be lower than in $\Gamma$, and thus passing from $\Gamma$ to $\widehat \Gamma$ should actually help. But, of course, in order to be able to pass percolation information from $\Gamma$ to $\widehat\Gamma$, we still need the tiles to satisfy Properties \iref{i.2} and  \iref{i.3}; in particular, Question~\ref{q.unique} needs to be resolved positively.

For a scale-invariant, or more generally, residually finite, amenable group, there is a simple method to construct  in any of its Cayley graphs a sequence of {\bf F{\o}lner monotilings}, i.e., the $n$th tiling has a single connected prototile $T_n$ such that these prototiles form a F{\o}lner sequence exhausting the graph, thus settling Property~\iref{i.2} for these groups. This was explained to us by G\'abor Elek, but as we learnt later, almost the same result was proved by Benjy Weiss in \cite{Benjy}: any residually finite amenable group and any solvable group has a F{\o}lner sequence $T_n$ such that the group can be tiled with each $T_n$ as a single prototile. Weiss did not make sure that the $T_n$ are connected, but this needs only a small trick. We give a proof below, basically due to Elek. An intriguing open problem in \cite{Benjy} is whether every group has monotilings with the prototiles exhausting the group. 

\bth\label{t.Elek}
Let $G$ be a residually finite amenable group, with a family of finite index subgroups $G_n$ with $\bigcap_{n\geq 1}G_n$ being finite. Then any finitely generated Cayley graph $\Gamma$ of $G$ has a sequence of tilings $\widehat\Gamma^n=\{ g T_n: g \in G_n\}$ such that the prototiles $T_n$ are connected and form a F{\o}lner exhaustion of $\Gamma$. In particular, for $G=\BS(1,m)$ or $\F\wr\Z$ or $Sol$, we can have $G_n\simeq G$.
\eth

\proof First assume for simplicity that $\bigcap_n G_n=\{1\}$. By taking $G_n^*:=\bigcap_{k\leq n} G_k$, we may assume that $G_n \supset G_{n+1}$ for all $n$. Consider the right coset tree $\T$ of $G$ w.r.t.~this subgroup chain. $G$ acts on $\T$ from the right, the $n$th level stabilizers are isomorphic to $G_n$, and the action is free on $\p\T$. So, the Schreier graph $\Gamma$ on $\p\T$ w.r.t.~some generating set $S$ of $G$ is in fact the Cayley graph $\Gamma(G,S)$, and the $n$th level Schreier graph $\Gamma_n$ with these generators is a finite factor-graph of $\Gamma$. Clearly, these $\Gamma_n$ converge locally to $\Gamma$. Since $\Gamma \twoheadrightarrow \Gamma_n$ is a topological covering map, there exists a 1-to-1 pre-image of the vertices of $\Gamma_n$ in $\Gamma$ inducing a connected subgraph $C_n$, and then $\{gC_n : g\in G_n\}$ is a tiling of $\Gamma$. However, it will not be true for any choice of $C_n$ that they converge random-locally to $\Gamma$.

The following version of the Ornstein-Weiss {\bf quasi-tiling lemma} \cite{OW} is proved in \cite{Elek:strong}: for any $\eps>0$, $N\in\N$, and a F{\o}lner exhaustion $F_n$ of an amenable Cayley graph $\Gamma$, there exist $\delta>0$, $L\in \N$, and a finite sub-collection of F{\o}lner sets $(F_{n_i})_{i=1}^s$ inside the ball $B_L(\Gamma)$ of radius $L$ of $\Gamma$, with $n_i > N$  $\forall i$, such that if $\Lambda$ is any finite graph with at least $1-\delta$ proportion of its vertices having an $L$-neighborhood isomorphic to $B_L(\Gamma)$, then $\Lambda$ can be $\eps$-quasi-tiled with translates of the $F_{n_i}$. In particular, in our sequence $(\Gamma_n)$ of finite graphs converging locally to $\Gamma$, for $n$ large enough, we can remove $O(\eps)$ proportion of the edges of $\Gamma_n$ such that each resulting component will be a subgraph of $B_L(\Gamma)$. (Thus the sequence $(\Gamma_n)$ is {\bf hyperfinite} in the sense of \cite{Elek:cost}.) Contract now each component into a single vertex, and using the edges between the components choose a spanning tree on the resulting graph. Then, along this spanning tree, we can lift each component from $\Gamma_n$ to $\Gamma$, keeping the pre-images connected in $\Gamma$ using the pre-images of the spanning tree edges. Thus we get a connected pre-image $T_n$ of $\Gamma_n$. The boundary edges $\p T_n$ are all pre-images of the edges between the components of $\Gamma_n$, and each such edge is covered by at most two edges of $\p T_n$, since it has two endpoints in $\Gamma_n$. Thus the boundary-to-volume ratio of $T_n$ in $\Gamma$ is at most $O(\eps)$. In summary, the $T_n$ are connected F{\o}lner sets, hence they converge random-locally to $\Gamma$, as desired.

If $\bigcap_{n} G_n=F$ is a finite group, then the Schreier graph  $\Gamma(G,\p\T,S)$ is a factor graph of $\Gamma(G,S)$, by an $|F|$-to-1 factor map $\pi$. The above argument can be run for $\Gamma(G,\p\T,S)$, then the tiles $T_n$ can be lifted by $\pi^{-1}$: any connected pre-image will form a F{\o}lner sequence in $\Gamma(G,S)$.
\qed
\medskip

Finally, what are the F{\o}lner sequences $F_n$ for which Question~\ref{q.unique} might have a positive answer, thus satisfying Property~\iref{i.3}? As pointed out e.g.~in \cite{antaretal}, for such questions it is important that there should exist some absolute constant $k$ such that the vertex boundary of $F_n$ is connected in the distance $k$ Rips complex, for all $n$. This is the case for any F{\o}lner sequence when $\Gamma$ is the Cayley graph of a finitely presented group. On the other hand, the usual F{\o}lner sequence of the lamplighter group $\Z_2 \wr \Z$ is not such, but it is not very difficult to augment each $F_n$ with some paths such that the resulting sequence $F_n^*$ already has this property  \cite{antaretal}.

As we mentioned in Subsection~\ref{ss.LL}, a natural Cayley graph of $\Z_2\wr \Z$ from the scale-invariance point of view is the Diestel-Leader graph $DL(2,2)$. Although it would be very interesting to do percolation renormalization on $DL(2,2)$, let us remark that the analogue of the $\Z^d$ half-space result is actually known already on $DL(2,2)$, using different methods \cite{ariel}. On the other hand, it is not known on any Cayley graph $\Gamma$ of $\Z_2\wr \Z$ that the unique infinite percolation cluster inherits the transience of simple random walk for all $p>p_c(\Gamma)$. The sequence $\widehat\Gamma^n$ of ``growing'' tiling graphs might be good enough e.g.~for this problem, provided the tiles satisfy Question~\ref{q.unique}. See \cite{BLS:pertu,ChPePe,repulsion} for more on the survival of random walk properties under percolation on groups.

\section{Scale-invariant tilings for the Heisenberg group}\label{s.Heis}

The main goal of this section is to construct a tiling sequence in the integer Heisenberg group satisfying Properties~\iref{i.1} and~\iref{i.2} of Section~\ref{s.perc}, i.e., to prove the following:

\bth\label{t.Heis} The discrete 3-dimensional Heisenberg group has Cayley graphs with strongly scale-invariant tilings (as deÞned in Question~\ref{q.tiling}). Moreover, the growing tiles form a F{\o}lner sequence.
\eth

We will give two explicit examples, both using the same somewhat general strategy, based on especially nice self-similar actions (see Definition~\ref{d.ss}) of the Heisenberg group that come from expanding endomorphisms of it (as defined in the third paragraph of the Introduction). For this, we have to start with several definitions and lemmas, culminating in the proof of the general Theorem \ref{t.everything} below. Most of the ideas are already contained in \cite{Nekra}, but the theorem itself is not. 
%\medskip

\subsection{The general strategy}

Given an {\bf expanding endomorphism} $\phi$  of the real
Heisenberg group $\mathcal{G}$ with the property that
$\phi(G)\subseteq G$ for the integer Heisenberg group $G$, and
$[G:\phi(G)]=t<\infty$, we can view it as a $t$-fold self-covering
of the compact Riemannian manifold $\mathcal{G}/G$. The inverse of
$\phi$ is an isomorphism $\psi: \phi(G) \longrightarrow G$; since $\phi(G)$ is of finite index,
$\psi$ is called a {\bf virtual endomorphism} of $G$. 
From this surjective virtual endomorphism, plus any transversal set
of coset representatives $\{g_0,\dots,g_{t-1}\}$ for $G/\phi(G)$, following \cite[Section
2.5.5]{Nekra}, one can get a self-similar action of $G$ on the
$t$-ary tree $\T=\{0,\dots,t-1\}^*$: for any $i\in\{0,\dots,t-1\}$, $w\in \T\cup \p\T$ and $g\in G$,
\be\label{e.action}
(iw)^g=jw^{\psi(g_j^{-1}gg_i)}\,, \quad\textrm{ where $j\in\{0,\dots,t-1\}$ is such that }g_j^{-1}gg_i \in \phi(G)\,.
\ee
Since $\phi$ is expanding, by
\cite[Theorem 6.1.3]{Nekra} we have that this self-similar action
of $G$ is {\bf contracting} in the following sense \cite[Section
2.11]{Nekra}: there exists a finite set $\NN \subset G$ such that
for every $g\in G$ there is a $k\in\N$ such that the restriction
$g|_{v}$ (as defined after Definition~\ref{d.ss}) is in $\NN$ for all vertices $v\in\T$ of depth at least
$k$. The minimal set $\NN$ with this property is called the {\bf
nucleus} of the self-similar action. It is easy to see that, because of the minimality of $\NN$, any restriction of a $g\in\NN$ is in the nucleus again, hence $\NN$ is a self-similar generating set of $\langle\NN\rangle$.

Given any contracting action by $G$ on the $t$-ary tree $\T$, one
can define the {\bf limit space} $\JJ_G$ as the quotient of the
set of left-infinite sequences $\{0,\dots,t-1\}^{-\N}$ by the
following {\bf asymptotic equivalence} relation: the sequences
$(\dots,x_{-1},x_{0})$ and $(\dots,y_{-1},y_{0})$ are equivalent
if{f} there exists a finite subset $K\subset G$ such that for all
$k\in\N$ there is some $g_k\in K$ with
$(x_{-k},\dots,x_{0})^{g_k}=(y_{-k},\dots,y_{0})$ with the
action of $G$ on $\T$, i.e., with $x_{-k}$ on the first level,
$(x_{-k},x_{-k+1})$ on the second level, etc. (In particular, this equivalence is very different from two rays in $\p \T=\{0,\dots, t-1\}^{\N}$ being in the same $G$-orbit.) A similar notion is the {\bf limit solenoid} $\SS_G$, the quotient of
$\{0,\dots,t-1\}^\Z$ by the equivalence relation that
$(\dots,x_{-1},x_0,x_1,\dots) \sim (\dots,y_{-1},y_0,y_1,\dots)$
if{f} there is a finite $K\subset G$ such that $\forall k\in\N$ 
$\exists g_k\in K$ with
$(x_{-k},x_{-k+1},\dots)^{g_k}=(y_{-k},y_{-k+1},\dots)$ in $\p\T$. Both on $\JJ_G$ and $\SS_G$, the topology is the image of the product topology under the equivalence quotient map.

We will look at the {\bf leaves} of $\SS_G$, 
corresponding to $G$-orbits $O\subset \p\T$: let $\LL_O$ be the
image of the set $\big\{\overline{x} \in \{0,\dots,t-1\}^\Z :
[\overline{x}]=(x_0,x_1,\dots)\in O \big\}$ in $\SS_G$. 
Note that if $\overline{x}$ and $\overline{y}$ are asymptotically equivalent, then, in particular, $[\overline{x}]$ and $[\overline{y}]$ are in the same $G$-orbit, hence different leaves are disjoint. For the topology on a leaf we do not take the restriction of the usual topology of the solenoid; rather, it is the image of the topology on $\{0,\dots,t-1\}^\Z$ that is product topology on the left tail but discrete on the right. The solenoid with this topology will be denoted by $\SS_G^\triangleleft$.

For each $w \in \p\T$, we define the {\bf tile}  $T_w\subseteq \SS_G^\triangleleft$ associated to $w$, 
the image under the quotient map of the set of sequences $\overline{x}$ with right tail
$[\overline{x}]=w$. For any
$w\in\p\T$, the set $\{T_{w^g} : g\in G \}$ is a covering of the
leaf $\LL_{O(w)} \subset \SS_G^\triangleleft$ corresponding to the $G$-orbit $O(w)$ of
$w$. However, the name ``tile'' is a bit misleading: because of the factorization, for two different
translates $g,h\in G$ the corresponding tiles do not have to be disjoint at all, in general.

\bl\label{l.tiles} For the tiles in the solenoid $\SS_G^\triangleleft$ of a contracting action $G$ on $\T$:\\
{\bf (i)} For $w\in\partial \T$, we have $T_{w^g} \cap T_{w^h} \not=\emptyset$ if and
only if $g^{-1}h$ is in the nucleus $\NN$.\\
{\bf(ii)} The tiles $T_w$ are connected if{f} the following graph is
connected. The vertices are $\{0, \ldots, t-1\}$, and $(i,j)$ is an edge if there exists an element $g\in\NN$ such that $(iw)^g=jw$ for all words $w$.\\
{\bf (iii)} If the {\bf open set condition} holds, i.e.,
$\forall g\in \NN$ $\exists v\in\T$ such that the restriction
$g|_v$ is the identity, then each tile is the closure of its
interior, and different tiles have disjoint interiors.
\el

\proof Part (i) is contained in \cite[Proposition 3.3.5]{Nekra}, 
but let us present here a self-contained proof to see how the definitions work, 
and since we will later use the argument of the ``if'' part again.

For the ``only if'' direction, assume that $\overline{x},\overline{y} \in \{0,\dots,t-1\}^\Z$  are asymptotically equivalent, with a finite set $K\subset G$ in the definition. Since the action is contracting, with nucleus $\NN$, for every $g \in K$ there is an $\ell_g$ such that $g|_v\in\NN$ for all $v\in\T$ of depth at least $\ell_g$. Taking $\ell:=\max\{\ell_g : g \in K\}<\infty$, we get that $(x_{-\ell},x_{-\ell+1},\dots)^{g_\ell}=(y_{-\ell},y_{-\ell+1},\dots)$ and $(x_0,x_1,\dots)^h=(y_0,y_1,\dots)$, where $h=g_\ell |_v \in \NN$ for $v=(x_{-\ell},x_{-\ell+1},\dots, x_1) \in \T$. Thus we are done. 

For the ``if'' part, $h_0:=g^{-1}h\in\NN$, hence, by the minimality of the nucleus, there is $h_1\in\NN$ and $x_{-1}\in\{0,\dots,t-1\}$ such that $h_1|_{x_{-1}}=h_0$, and $(x_{-1}w^g)^{h_1}=y_{-1}w^h$ for some $y_{-1}\in\{0,\dots,t-1\}$. Then there is $h_2\in\NN$ and $x_{-2}$ such that $h_2|_{x_{-2}}=h_1$, and so on.  This way, we get words $\overline{x}=(\dots,x_{-2},x_{-1})w^g$ and $\overline{y}=(\dots,y_{-2},y_{-1})w^h$ such that $\big((x_{-k},\dots,x_{-1})w^g\big)^{h_k}=(y_{-k},\dots,y_{-1})w^h$ for $h_k\in\NN$, $\forall k\in \N$. Since  $\NN$ is finite, this shows that $\overline{x}$ and $\overline{y}$ are asymptotically equivalent, so they map to one element in $T_{w^g}\cap T_{w^h}$.

Part (ii) is proved in~\cite[Proposition 3.3.10]{Nekra}. Note that a tile $T_w$ is the union of the $t$ sub-tiles $\big\{ (\dots,x_{-3},x_{-2}, i) w \big\}$, $i=0,\dots,t-1$, so the connectedness of $T_w$ has to do with the intersections between these sub-tiles and the sub-tiles of those, and so on, hence part (i) is of relevance here. 

Part (iii) is \cite[Proposition 3.3.7]{Nekra}. \qed

\bl\label{l.leaf}
Points in different leaves cannot be connected in $\SS_G^\triangleleft$, i.e., every leaf is a union of path-connected components of the solenoid.
\el

\proof For $w_1\not=w_2 \in \p\T$, the sets $\big\{ (\dots,x_{-2},x_{-1}) w_i \big\}$ that map to $T_{w_i} \subset \SS_G^\triangleleft$ are disjoint clopen sets before the factorization, and, by Lemma~\ref{l.tiles}~(i), they can have a common point after the factorization only if $w_1$ and $w_2$ are in the same $G$-orbit. \qed
\medskip

For a self-similar action of $G$ on $\T$, the wreath product decomposition~(\ref{e.treathG}) defines a homomorphism from each first level stabilizer $\St_G(i)$, $i=0,1,\dots,t-1$, to $G$. If the action is transitive on the first level, then the $\St_G(i)$ are all conjugate subgroups of $G$ with index $t$, so these endomorphisms $\St_G(i) \longrightarrow G$ are also conjugates of each other, and they are virtual endomorphisms of $G$.

A self-similar action of $G$ is called {\bf self-replicating} ({\bf recurrent} in~\cite{Nekra}) if it is transitive on the first level, and the associated virtual endomorphisms of $G$ are onto. When the action is constructed from an expanding homomorphism $\phi$ of $G$, then it is clear from (\ref{e.action}) that the virtual endomorphisms of $G$ associated to the self-similar action are conjugates of the virtual endomorphism $\psi=\phi^{-1}$, and the action is self-replating.

\bl\label{l.selfrep}
Assume that $G$ has a contracting self-replicating action on $\T$. Then {\bf (i)} the nucleus $\NN$ generates $G$; {\bf (ii)} the limit space $\JJ_G$ is path-connected and locally path-connected; {\bf (iii)} the leaves $\LL_O$ of the limit solenoid $\SS_G^\triangleleft$ are exactly its path-connected components.
\el

\proof Part (i) is easy from the definitions of contracting and self-replicating; see \cite[Proposition 2.11.3]{Nekra} for details. Parts (ii) and (iii) are Theorem 3.6.3 and Proposition 5.7.9 of \cite{Nekra}, respectively, plus the fact that for locally compact metrizable spaces connectedness plus local connectedness implies the path-connected versions, see \cite[Corollary 3.5.3]{Nekra}. Note that they should not be surprising in light of Lemmas~\ref{l.tiles}~(i), (ii) and~\ref{l.leaf} together with $\langle \NN \rangle=G$. \qed
\medskip

Summarizing: for a contracting self-replicating action satisfying the conditions of Lemma~\ref{l.tiles} (connected tiles and the open set condition), on the connected leaf $\LL_O$ we obtain a tiling (in the usual geometric sense) by the connected tiles $T_w$, $w\in O$. By Lemma~\ref{l.tiles}~(i),
the adjacency graph of this tiling is the Schreier graph $\Gamma(G,O,\NN)$ of the action
of $G$ on the orbit $O\subset \p\T$, with generators $\NN$. 
If the $G$-stabilizer of $w\in O$ is trivial, then the
Schreier graph is in fact the Cayley graph
$\Gamma(G,\NN)$, which is connected since $\langle \NN \rangle=G$.

Now consider the shift map $\shift$ that moves the origin to the left in
$\{0,\dots,t-1\}^{\Z}$, or deletes the last letter in
$\{0,\dots,t-1\}^{-\N}$ (hence, it is $t$-to-1). In both cases, $\shift$ preserves the
asymptotic equivalence relation, and thus we get the dynamical
systems $(\JJ_G,\shift)$,  $(\SS_G,\shift)$ and $(\SS_G^\triangleleft,\shift)$.
When the contracting action is obtained from an expanding
endomorphism $\phi$, the above quoted  \cite[Theorem 6.1.3]{Nekra}
also says that $(\JJ_G,\shift)$ is topologically conjugate to
$(\mathcal{G}/G,\phi)$, moreover, the {\bf iterated monodromy group} $\mathsf{IMG}(\phi)$
of the $t$-fold self-covering $\phi: \mathcal{G}/G \longrightarrow
\mathcal{G}/G$ is exactly $G=\pi_1(\mathcal{G}/G)$, with a self-similar action that is basically~(\ref{e.action}).  
(So, $\mathsf{IMG}$ and the limit space constructions are inverses of each other.) 
Furthermore, as it will become clear in the next paragraph, $(\SS_G^\triangleleft,\shift)$ resembles $(\mathcal{G},\phi)$ in some sense, but they are certainly not topologically conjugate, since $\SS_G^\triangleleft$ is highly disconnected, and we cannot restrict $\shift$ to a given leaf $\LL_O$, either, since a shift typically leaves the leaf. We will use the dynamical system $(\SS_G^\triangleleft,\shift)$, not $(\mathcal{G},\phi)$ or $(\mathcal{G}/G,\phi)$, but the latter two are certainly easier to visualize, so the Reader is encouraged to keep them in mind.

The shift map $\shift$ gives a grouping of the tiles in
$\SS_G^\triangleleft$:  for any $w\in\p\T$, we have
$$\shift(T_w)=\bigcup_{0 \leq i < t} T_{iw}\,.$$
One hopes that this gives a scaling of tiling graphs, as in
Question~\ref{q.tiling}. However, the points $w$ and $iw$ might not be in the same $G$-orbit, hence $T_w$ and $T_{iw}$ could be in different leaves. If the action is self-replicating, then it is easy to see that at least the $iw$ are all in the same orbit, and if also the ``connected tiles'' condition of Lemma~\ref{l.tiles}~(ii) is satisfied, together they form a connected set in the tiling graph which is just the Schreier graph $\Gamma(G,O(0w),\NN)$. So, after the grouping  we have the connected new tiles $\shift(T_w)$. By the minimality of $\NN$, as in Lemma~\ref{l.tiles}~(i), given $w_1,w_2\in\p\T$, there exists $n\in\NN$ with $w_1^n=w_2$ if and only if there exist $i,j\in\{0,\dots,t-1\}$ and $m\in\NN$ such that $(iw_1)^m=jw_2$. Therefore, by Lemma~\ref{l.tiles}~(i), the new tiles $\shift(T_{w^g})$ and $\shift(T_{w^h})$ are neighbors in $\LL_{O(0w)}$ if{f} $T_{w^g}$ and $T_{w^h}$ are neighbors in $\LL_{O(w)}$. That is, the new tiling graph will be isomorphic to $\Gamma(G,O(w),\NN)$. However, this graph might be different from $\Gamma(G,O(0w),\NN)$. So, we would also like that $G$ acts
freely on $\p\T$: this way, the Schreier graphs
$\Gamma(G,O(w),\NN)$ will be isomorphic to the Cayley graph
$\Gamma(G,\NN)$ for all $w\in\p\T$. 

We are now ready to state our general theorem.

\bth\label{t.everything}
If $G$ has a contracting self-replicating action on the $t$-ary tree $\T$ with the  ``connected tiles''  and ``open set'' conditions, and the action on $\p\T$ is free, then $G$ has a Cayley graph with a strongly scale-invariant tiling sequence, in which the growing tiles form a F{\o}lner sequence. 
\eth

\proof We have just said how one grouping of tiles is done, but how do we get a sequence of tilings in the same Cayley graph, with tiles exhausting it? Given a ray $w=w_0w_1\dots \in\p\T$, consider the grouping
$$
\shift^{n}(T_{w_nw_{n+1}\dots})=\bigcup \big\{ T_{i_0\dots i_{n-1}w_nw_{n+1}\dots} \, : \, i_0,\dots,i_{n-1} \in \{0,\dots,t-1\} \big\}\,.
$$
As discussed above, the tiles on the right hand side are all in the same leaf $\LL_{O(w)}$. So, this grouping for $n=1,2,\dots$ gives larger and larger connected tiles in the Schreier graph $\Gamma(G,O(w),\NN)$, and if the action of $G$ is free on $\p\T$, then each tiling graph is isomorphic to the Cayley graph $\Gamma(G,\NN)$. But we still have to pick $w$ cleverly to have these tiles exhaust the graph $\Gamma(G,O(w),\NN)$, i.e., we want that for all $w'\in O(w)$ there is $n_0\in\N_+$ such that $w'_n=w_n$ for all $n\geq n_0$.

\bl\label{l.random}
Assume that a contracting action on the $t$-ary tree with nucleus $\NN$ satisfies the open set condition. Consider an i.i.d.~random sequence $(\xi_n)_{n=1}^\infty$ with $\xi_n\in\{0,1,\dots,t-1\}$ chosen uniformly. Then there are some $c,C>0$ such that $\Pr\big[\exists g\in \NN \textrm{ with } g|_{\xi_1\dots \xi_n}\not= \mathrm{id} \big] \leq C\exp(-cn)$ for all $n\in \N_+$. 
\el

\proof Consider the Moore diagram of the nucleus $\NN$. The open set condition means that $\mathrm{id}\in\NN$ is accessible from each element $g\in\NN$ with a finite directed path in the diagram; let the maximum length of these paths over the starting points $g$ be $\ell$. Note that $g|_{\xi_1\dots \xi_n}$ is nothing else but the state in $\NN$ after $n$ uniform random steps in the Moore diagram starting from $g$. Starting from anywhere, the probability of reaching id in $\ell$ steps is at least $t^{-\ell}$, hence the probability of not reaching id in $n$ steps is at most $(1-t^{-\ell})^{\lfloor n/\ell\rfloor}$, which is at most $C\exp(-cn)$ for some $C,c>0$. So, the probability that $g|_{\xi_1\dots \xi_n}\not=\mathrm{id}$ for some $g\in \NN$ is at most $|\NN|C\exp(-cn)$, and we are done.\qed
\medskip

Since every $g\in G$ is a product of finitely many elements from $\NN$, the lemma implies that almost every random sequence $w=(\xi_n)_{n=1}^\infty$ is such that for all $g\in G$ there is an $n_0\in \N_+$ with $(w^g)_n = w_n$ for all $n\geq n_0$, as required for the tiles to exhaust the graph $\Gamma(G,O(w),\NN)$.

The last thing to check is that the tiles form a F{\o}lner sequence. The tile $\shift^{n}(T_{w_nw_{n+1}\dots})$ has $t^n$  of the original  small tiles, hence $t^n$ vertices of $\Gamma(G,O(w),\NN)$. By Lemma~\ref{l.tiles}~(i), such a small tile  $T_{i_0\dots i_{n-1}w_nw_{n+1}\dots}$ can be on the boundary of the large tile only if there is a $g\in \NN$ with $g|_{i_0\dots i_{n-1}} \not= \mathrm{id}$. By Lemma~\ref{l.random}, the proportion of these tiles among all the $t^n$ is at most $C\exp(-cn)$, hence the growing tiles indeed form a F{\o}lner sequence, with polynomially small boundary-to-volume ratio. \qed
\medskip

\noindent{\bf Remark.} Both the existence of an exhausting tile sequence and the F{\o}lner condition relied on Lemma~\ref{l.random}. There is a general connection between these issues: given the recursively defined tiles $T_n$ in a transitive graph $\Gamma$ forming a F{\o}lner sequence, we can translate the tiles to get an exhausting sequence. The trick is that  for any $n$ there is $N=N(n)$ such that $T_N$ contains a copy of $T_n$ whose boundary edges are all inside $T_N$, as soon as the edge boundary to volume ratio of $T_N$ is less than $1/|T_n|$. So, fix $x\in \Gamma$, a tile $T_{n_0}\ni x$, then let $n_{k+1}=N(n_k)$ for $k\in\N$, and translate $T_{n_k}$ so that it contains $T_{n_{k-1}}$ in its interior. Then $T_{n_k}$ will contain the ball $B_k(x)$, hence they exhaust $\Gamma$.

\subsection{The proof of Theorem~\ref{t.Heis}} 
%\proof[Proof of Theorem~\ref{t.Heis}] 

Starting from an expanding endomorphism $\phi$ ensures that the action (\ref{e.action}) is self-replicating. So, to apply Theorem~\ref{t.everything}, in the specific examples we will need that the action is free on $\p\T$ and satisfies the ``open set'' and ``connected tiles'' conditions.
\medskip

As the first example, consider the expanding Heisenberg group
endomorphism
$$
\phi:\  \left(\begin{array}{ccc} 1
& a & b \\ 0 & 1 & c\\ 0 & 0 & 1\end{array}\right)
\mapsto
\left(\begin{array}{ccc} 1
& 2c & -2b+2ac\\ 0 & 1 & a\\ 0 & 0 & 1\end{array}\right),
$$
with $[G:\phi(G)]=4$. Denote
\[A=\left(\begin{array}{ccc}1 & 1 & 0\\ 0 & 1 & 0\\ 0 & 0 &
1\end{array}\right),\quad C=\left(\begin{array}{ccc}1 & 0 & 0\\ 0
& 1 & 1\\ 0 & 0 & 1\end{array}\right),\quad
B=\left(\begin{array}{ccc}1 & 0 & 1\\ 0 & 1 & 0\\ 0 & 0 &
1\end{array}\right).\] 
The inverse of $\phi$ is the virtual
endomorphism that maps $A^2$ to $C$, $C$ to $A$ and $B^2$ to
$B^{-1}$. Using (\ref{e.action}), we get the
following self-similar action on the 4-ary tree, using the notation implied by (\ref{e.treathG}):
\[A=(1, C, 1, C)(01)(23)\]
\[C=(A, CAC^{-1}, A, A)(13).\]
The commutator $B=[A, C]=A^{-1}C^{-1}AC$ decomposes then as
\[B=(C^{-1}A^{-1}CA, A^{-1}CAC^{-1}, 1, 1)(02)(13)=(B^{-1}, CB^{-1}C^{-1}, 1, 1)(02)(13).\]
It follows that
\[[A, B]=(1, 1, BCB^{-1}C^{-1}),\]
\[[C, B]=(1,  C\cdot B^{-1}A^{-1}BA\cdot C^{-1}, A^{-1}BAB^{-1}, 1),\]
which gives an inductive proof of the relations $[A, B]=1$ and
$[C, B]=1$.

Furthermore, we have
\[CAC^{-1}=(1, C, ACA^{-1}C^{-1}, ACA^{-1})(03)(12)=
(1, C, B, ACA^{-1})(03)(12)\]
\[B=(B^{-1}, B^{-1}, 1, 1)(02)(13)\]
\[ACA^{-1}=(CAC^{-1}, CAC^{-1}, A, CAC^{-1})(02).\]
Therefore, $\{1, A, C, CAC^{-1}, B, B^{-1} ACA^{-1}\}$ is
a self-similar generating set, and every element of this set has a
trivial restriction, which implies that the group satisfies the open
set condition.

It is easy to see that no
subgroup of $G$ is invariant under the action of the virtual endomorphism,
hence the above action of the Heisenberg group on the tree is faithful
and free on the boundary.

The nucleus can be obtained by computer (using the GAP
package~\cite{MS:gap} developed for this purpose): it consists of
the trivial element, the elements
\[A, C, CA, AC, AC^{-1}, B,  BA, BC, BA^{-1}, BC^{-1},
BCA^{-1}, BAC^{-1},\] and their inverses; 25 elements in total.
Note that  $A(0w)=1w$, $A(2w)=3w$, $B(2w)=0w$ and
$B(3w)=1w$, which implies that the tiles are connected.
Therefore, the Cayley graph of $G$ generated by these elements has a strongly scale-invariant tiling with F{\o}lner tiles.
\medskip

Our second example starts with the expanding endomorphism that Gelbrich used in
\cite{GelbHeis} to obtain a periodic self-similar tiling of the real Heisenberg group:
\[\left(\begin{array}{ccc}1 & a & b\\ 0 & 1 & c\\ 0 & 0 &
1\end{array}\right)\mapsto\left(\begin{array}{ccc} 1 & a+c & 2b-ac+\frac{c^2-c-a^2+a}2\\
0 & 1 & c-a\\ 0 & 0 & 1\end{array}\right),\] again with
$[G:\phi(G)]=4$.

The associated self-similar action is given by
\[A=(1, C^{-1}A, 1, C^{-1}A)(01)(23),\]
\[C=(C, A, C, C^{-1}AC)(0123) =  (C, A, C, AB)(0123),\]
where we denote again $B=[A, C]=A^{-1}C^{-1}AC$. 
The restictions of $A$ and $C$ are decomposed as
\[C^{-1}A=(C^{-1}, C^{-1}, BC^{-1}, C^{-1})(02)\]
\[AB=(1, C^{-1}A, B, C^{-2}AC)(03)(12)\]
\[BC^{-1}=(A^{-1}, C^{-1}, A^{-1}, BC^{-1})(0123)\]
\[B=(1, 1, B, B)(02)(13)\]
\[C^{-2}AC=(BC^{-1}, C^{-1}, BC^{-1},
BC^{-1})(13).\]

We see that the set \[\{1, A, B, C,
AB, BC^{-1}, C^{-1}A, C^{-2}AC\}^{\pm 1}\] is self-similar and
that the identity is accessible from each of these elements as a restriction.

The nucleus consists of the trivial element, the elements
\[A, B, C, AB, BC, BC^{-1}, C^{-1}A, C^{-2}AC\]
and their inverses (17 elements in total). The tiles are
connected, since $A(0w)=1w$, $A(2w)=3w$, $B(0w)=2w$ and $B(1w)=3w$
for all words $w$. Hence all conditions of Theorem~\ref{t.everything} are satisfied.\qed
\medskip

\noindent{\bf Acknowledgments.} We are grateful to Itai Benjamini,
G\'abor Elek, David Fisher, Mark Sapir and \'Ad\'am Tim\'ar for helpful conversations and remarks on the manuscript, 
to Douglas Lind and Mark Sapir for bringing \cite{Benjy} and \cite{Bridsonetal} to our attention, and to the referees for their suggestions.

The collaboration between the authors began at an American Institute
of Mathematics workshop in Palo Alto in October 2007. The research of V.N.~is
supported by NSF grant DMS-0605019. During most of this work,
G.P.~was a postdoc at Microsoft Research, Redmond, and at MSRI, Berkeley, and 
was partially supported by the Hungarian OTKA, grant T049398; presently supported by an NSERC Discovery Grant.


\begin{thebibliography}{BartGN03}

\bibitem[AlLy07]{urn} D. Aldous and R. Lyons.
Processes on unimodular random networks.
{\it Electron. J. Probab.} {\bf 12} (2007), Paper 54, 1454--1508. 

\bibitem[AlBS04]{ABS} N. Alon, I. Benjamini and A. Stacey.
Percolation on finite graphs and isoperimetric inequalities. {\it Ann. Probab.} {\bf 33} (2004), 1727--1745.

\bibitem[AnP96]{chemical} P. Antal and A. Pisztora.
On the chemical distance in supercritical Bernoulli percolation.
{\it Ann. Probab.} {\bf 24} (1996), 1036--1048.

\bibitem[BanST08]{antaretal} A. Bandyopadhyay, J. Steif and \'A. Tim\'ar.
On the cluster size distribution for percolation on some general graphs.
{\it Preprint,} {\tt arXiv:0805.3620 [math.PR]}.

\bibitem[BarsGN91]{half} D. Barsky, G. Grimmett and C. Newman. Percolation in halfspaces:
equality of critical densities and continuity of the percolation probability.
{\it Probab. Theory Related Fields} {\bf 90} (1991), no. 1, 111--148.

\bibitem[BartGN03]{fractal} L. Bartholdi, R. Grigorchuk and V. Nekrashevych.
From fractal groups to fractal sets. {\it Fractals in Graz (P. Grabner and W. Woess, eds.)}, Trends in Mathematics, Birkh\"auser Verlag, Basel, 2003, pp. 25--118. {\tt arXiv:math.GR/0202001v4}

\bibitem[BartG\v{S}03]{branch} L. Bartholdi, R.~I. Grigorchuk and Z. \v{S}uni\'k.
Branch groups. {\it Handbook of Algebra} vol 3, pp. 989--1112. North-Holland, Amsterdam, 2003.
{\tt arXiv:math.GR/0510294v2}

\bibitem[Bart\v{S}06]{solvauto} L. Bartholdi and Z. \v Suni\'k. Some solvable automaton groups.
In {\it Topological and asympotic aspects of group theory}, {\it Contemp. Math.} {\bf 394} (2006), pp. 11--30.

\bibitem[BauS62]{BS} G. Baumslag and D. Solitar.
Some two-generator one-relator non-Hopfian groups.
{\it Bull. Amer. Math. Soc.} {\bf 68} (1962), 199--201.

\bibitem[BekV97]{BV} M. Bekka and A. Valette.
Group cohomology, harmonic functions and the first $L^2$-Betti number.
{\it Potential Anal.} {\bf 6} (1997), no. 4, 313--326.

\bibitem[Bel03]{Beleg} I. Belegradek.
On co-Hopfian nilpotent groups.
{\it Bull. London Math. Soc.} {\bf 35} (2003) 805--811.

\bibitem[Ben06]{Itai} I. Benjamini. Post on personal website, {\tt http://www.wisdom.weizmann.ac.il/$\sim$itai}, 2006.

\bibitem[BeLS99]{BLS:pertu} I. Benjamini, R. Lyons and O. Schramm.
Percolation perturbations in potential theory and random
walks, {\it In: Random walks and discrete potential theory (Cortona,
1997)}, Sympos. Math. XXXIX, M. Picardello and W. Woess (eds.),
Cambridge Univ. Press, Cambridge, 1999, pp. 56--84. {\tt arXiv:math.PR/9804010}

\bibitem[BenS96]{BS:beyond} I. Benjamini and O. Schramm. Percolation beyond
$\Z^d$, many questions and a few answers, {\it Elect. Commun. Probab.}
{\bf 1} (1996), 71--82.

\bibitem[BenS01]{BS:limit} I. Benjamini and O. Schramm.
Recurrence of distributional limits of finite planar graphs.
{\it Electron. J. Probab.} {\bf 6}, no. 23, 13 pp. (electronic).

\bibitem[BriHM07]{Bridsonetal} M. Bridson, A. Hinkkanen and G. Martin.
Quasiregular self-mappings of manifolds and word hyperbolic groups.
{\it Compos. Math.} {\bf 143} (2007), no. 6, 1613--1622.

\bibitem[BruSi98]{affine} A.~M. Brunner and S. Sidki.
The generation of $GL(n,\Z)$ by finite state automata.
{\it Internat. J. Algebra Comput.} {\bf 8 }(1998), no. 1, 127--139.

\bibitem[Bry07]{Brydges} D. Brydges. {\it Lectures on the Renormalisation Group.}
Park City Mathematics Institute lecture notes, 2007.
{\tt http://www.math.ubc.ca/$\sim$db5d/Seminars/}

\bibitem[ChPP04]{ChPePe} D. Chen and Y. Peres, with an appendix by G. Pete.
Anchored expansion, percolation and speed. {\it Ann. Probab.} {\bf 32}
(2004), 2978--2995. {\tt arXiv:math.PR/0303321}

\bibitem[DiL01]{DL} R. Diestel and I. Leader.
A conjecture concerning a limit of non-{C}ayley graphs.
{\it J. Alg. Combin.} {\bf 14} (2001), 17--25.

\bibitem[DrS08]{DruSap} C. Dru\c{t}u and M. Sapir.
Groups acting on tree-graded spaces and splittings of relatively hyperbolic groups.
{\it Adv. Math.} {\bf 217} (2008), no. 3, 1313--1367. {\tt arXiv:math.GR/0601305}.

\bibitem[Ele06]{Elek:strong} G. Elek. The strong approximation conjecture holds for amenable groups.
{\it J. Funct. Anal.} {\bf 239} (2006), no. 1, 345--355. {\tt arXiv:math.FA/0511655v1}

\bibitem[Ele07]{Elek:cost} G. Elek. The combinatorial cost.
{\it Enseign. Math.} {\bf 53} (2007), 225--235. {\tt arXiv:math.GR/0608474}

\bibitem[EFW07]{EFW} A. Eskin, D. Fisher and K.Whyte. Quasi-isometries and rigidity of solvable groups. {\it Pure and Applied Mathematics Quarterly} {\bf 3} (2007) 927--947.

\bibitem[Far81]{Farkas} D.~R. Farkas.
Crystallographic groups and their mathematics.
{\it Rocky Mountain J. Math.} {\bf 11} (1981), no.~4, 511--551.

\bibitem[Fra70]{Franks} J. Franks.
Anosov diffeomorphisms. {\it Global Analysis (Proc.~Sympos.~Pure Math., Vol.~XIV, Berkeley, Calif., 1968)} pp. 61--93. Amer. Math. Soc., Providence, R.I., 1970.

\bibitem[Gel94]{GelbHeis} G. Gelbrich.
Self-similar periodic tilings on the Heisenberg group.
{\it J. Lie Theory}  {\bf 4} (1994), no.~1, 31--37.

\bibitem[Gel95]{Gelb} G.~Gelbrich.
Self-similar tilings and expanding homomorphisms of groups.
{\it Arch. Math. (Basel)\/} {\bf 65} (1995), no.~6, 481--491.

\bibitem[Gra08]{binomial} A. Granville. Arithmetic properties of binomial coefficients.
Web survey. {\tt http://www.dms.umontreal.ca/ $\sim$andrew/Binomial/}

\bibitem[GrNS00]{GNS} R.~I. Grigorchuk, V.~V. Nekrashevich and V.~I. Sushchanskii.
Automata, dynamical systems and groups.
{\it Proceedings of the Steklov Institute of Mathematics} {\bf 231} (2000), 128--203.

\bibitem[Gr\.Z01]{GZ} R. Grigorchuk and A. \.Zuk. The lamplighter group as a group generated by a 2-state automaton and its spectrum. {\it Geom. Dedicata} {\bf 87} (2001), no.~1-3, 209--244.

\bibitem[Grim99]{Grimm}
G. Grimmett. {\it Percolation. Second edition}. Grundlehren der
Mathematischen Wissenschaften, 321. Springer-Verlag, Berlin, 1999.

\bibitem[GrimM90]{GrM}
G. Grimmett and J. Marstrand. The supercritical phase of percolation is well-behaved. {\it Proc. Roy. Soc. London Ser. A} {\bf 430} (1990), 439--457.

\bibitem[Gro81]{Gromov} M. Gromov.
Groups of polynomial growth and expanding maps. Appendix by J. Tits.
{\it Publ. Math. I.H.E.S.} {\bf 53} (1981), 53--73, 74--78.

\bibitem[dlH00]{dlH} P. de la Harpe.
{\it Topics in geometric group theory.}
Chicago Lectures in Mathematics, University of Chicago Press, 2000.

\bibitem[KaV83]{KV} V.~A. Kaimanovich and A.~M. Vershik.
Random walks on discrete groups: boundary and entropy.
{\it Ann. Probab.} {\bf 11} (1983), 457--490.

\bibitem[Kap08]{Kap:GGT} M. Kapovich: {\it Lectures on Geometric Group Theory}. 
{\tt http://math.ucdavis.edu/$\sim$kapovich/EPR/ggt.pdf}

\bibitem[Lyo00]{Ly:survey} R. Lyons. Phase transitions on nonamenable graphs.
{\it J. Math. Phys.} {\bf 41} (2000), 1099--1126. 

\bibitem[LyPP96]{LPP} R. Lyons, R. Pemantle and Y. Peres. Random walks on the lamplighter group. {\it Ann. Probab.} {\bf 24} (1996) no.~4, 1993--2006.

\bibitem[LyP08]{LPbook} R. Lyons, with Y. Peres. {\it Probability on trees and networks.} Book in preparation, present version is at {\tt http://mypage.iu.edu/$\sim$rdlyons}.

\bibitem[MS08]{MS:gap} Y. Muntyan and D. Savchuk. {\it \verb!AutomGrp! --- \verb!GAP! package
for computations in self-similar groups and semigroups, Version
1.1.2}, 2008. (available at
\verb!http://finautom.sourceforge.net!).

\bibitem[Nek05]{Nekra} V.~V. Nekrashevych. {\it Self-similar groups}.
Mathematical Surveys and Monographs, Vol.~117, Amer. Math. Soc., Providence, RI, 2005.

\bibitem[OrW87]{OW} D. S. Ornstein and B. Weiss.
Entropy and isomorphism theorems for actions of amenable groups.
{\it  J. Anal. Math.} {\bf 48} (1987), 1--141.

\bibitem[Pas93]{Paschke} W.~L. Paschke.
A numerical invariant for finitely generated groups via actions on graphs.
{\it Math. Scand.} {\bf 72} (1993), no.~1, 148--160.

\bibitem[PPS06]{ariel} Y. Peres, G. Pete and A. Scolnicov.
Critical percolation on certain non-unimodular graphs.
{\it New York J. of Math.} {\bf 12} (2006), 1--18. {\tt arXiv:math.PR/0501532v3}

\bibitem[Pes08]{Pestov} V. G. Pestov.
Hyperlinear and sofic groups: a brief guide.
{\it  Bulletin of Symbolic Logic}, to appear. {\tt arXiv:0804.3968v8 [math.GR]}

\bibitem[Pet08]{repulsion} G. Pete.
A note on percolation on $\Z^d$: isoperimetric profile via exponential cluster repulsion.
{\it Elect. Comm. Probab.} {\bf 13} (2008), 377--392. {\tt arXiv:math.PR/0702474v4}

\bibitem[Sap07]{Sapir} M.~V. Sapir. Some group theory problems. {\it Internat. J. of Algebra
and Comput.} {\bf 17} (2007), 1189--1214. {\tt arXiv:0704.2899v1 [math.GR]}

\bibitem[Sel97]{Sela} Z. Sela.
Structure and rigidity in (Gromov) hyperbolic groups and discrete groups in rank 1 Lie groups, II.
{\it Geom. Funct. Anal.} {\bf 7} (1997), no.~3, 561--593.

\bibitem[SiSt05]{SiSt} P.~V. Silva and B. Steinberg.
On a class of automata groups generalizing lamplighter groups.
{\it Internat. J. Algebra and Comput.} {\bf 15} (2005), 1213--1235.

\bibitem[SW91]{SW} P. Soardi and W. Woess.
Uniqueness of currents in infinite resistive networks. 
{\it Discrete Appl. Math.}  {\bf 31} (1991), no.~1, 37--49.

\bibitem[Wei01]{Benjy} B. Weiss.
Monotileable amenable groups. {\it In: Topology, ergodic theory, real algebraic geometry}, 257--262,
Amer. Math. Soc. Transl. Ser. 2, 202, Amer. Math. Soc., Providence, RI, 2001.

\bibitem[Wil98]{profinite} J.~S. Wilson.
{\it ProÞnite groups.} 
London Mathematical Society Monographs. New Series, 19. 
The Clarendon Press, Oxford University Press, New York, 1998.

\bibitem[Woe05]{W:DL} W. Woess. Lamplighters, Diestel-Leader graphs, random walks, and harmonic functions. {\it Combin. Probab. \& Comput.} {\bf 14} (2005), 415--433.


\end{thebibliography}
\end{document}